%% file: noteBWconjecture.tex
\documentclass[a4paper]{amsart}
\pdfoutput=1
\usepackage[T2A,T1]{fontenc}
\usepackage[utf8]{inputenc}
\usepackage{lmodern,fixltx2e}%
\usepackage{etex}%
\usepackage[german,russian,british]{babel}
\usepackage{hyphxcpt}
\usepackage{amsmath}
\usepackage{amssymb}
\usepackage{amsfonts}
\usepackage{enumerate}
\usepackage{mathoperators}
\usepackage[algoruled,linesnumbered]{algorithm2e}
\usepackage[uglyemptyset,hyperref]{usefulthings}
\usepackage{PolInv}
\usepackage{bm}
\usepackage[mathscr]{eucal}
\DeclareMathOperator{\ppc}{PPC\#}%
\DeclareMathOperator{\cdeg}{cdeg}%
\newcommand{\ncent}[2][n]{\ensuremath{{#2^{*(#1)}}}}%
\newcommand{\centn}[2][n]{\ensuremath{{#2^{(#1)*}}}}%
\newcommand{\nbicent}[2][n]{#2^{**(#1)}}%
\newcommand{\bicentn}[2][n]{#2^{(#1)**}}%
\newcommand{\cent}[1]{\ensuremath{{#1^{*}}}}%
\newcommand{\bicent}[1]{#1^{**}}%
\def\graph{}%
\DeclareMathOperator{\fix}{fix}%
\DeclareMathOperator{\pr}{pr}%
\renewcommand{\graph}[1]{\ensuremath{#1^{\bullet}}}%
\newcommand{\downdiag}{\mathord{\swarrow}}
\newcommand{\updiag}{\mathord{\nearrow}}
\newcommand*{\asc}{\mathord{\Uparrow}}
\newcommand*{\desc}{\mathord{\Downarrow}}
\newcommand*{\Sqre}{\mathord{\boxdot}}
\newcommand{\bfa}[1]{{\bm#1}}%
\newtheorem{lemma}{Lemma}
\newtheorem{proposition}[lemma]{Proposition}
\newtheorem{theorem}[lemma]{Theorem}
\newtheorem{corollary}[lemma]{Corollary}
\newtheorem{remark}[lemma]{Remark}
\theoremstyle{definition}
\newtheorem{algo}[lemma]{Algorithm}
\newtheorem{example}[lemma]{Example}
\author{Mike Behrisch}
\title{A note on the Burris-Willard conjecture}
\address{In\-sti\-tute of Discrete Mathematics and Geometry\\
         Tech\-ni\-sche Uni\-ver\-si\-t\"{a}t Wien\\
         Wied\-ner Haupt\-str.~8--10/E104\\
         \mbox{A-1040} Vienna\\
         Austria}%
\email{\href{mailto:behrisch@logic.at}{behrisch@logic.at}}%
\urladdr{\url{https://orcid.org/0000-0003-0050-8085}}%
\begin{document}
\begin{abstract}
Based on results by Dani\v{l}\v{c}enko, in 1987 Burris and Willard have
conjectured that on any \nbdd{k}element domain
where $k\geq 3$ it is possible to bicentrically generate every
centraliser clone from its \nbdd{k}ary part.

Later, for every $k\geq 3$, Snow constructed algebras with a
\nbdd{k}element carrier set where the minimum arity of the clone of term
operations from which the bicentraliser can be generated is at least
$(k-1)^2$, which is larger than~$k$ for $k\geq 3$.

We prove that Snow's examples do not violate the
Burris-Willard conjecture nor invalidate the results by
Dani\v{l}\v{c}enko on which the latter is based.
We also complement our results with some computational evidence for
$k=3$, obtained by an algorithm to compute a primitive positive
definition for a relation in a finitely generated relational clone
over a finite set.
\end{abstract}

\maketitle

\section{Introduction}
Centraliser clones are collections of homomorphisms of finite powers of
algebras into themselves. That is, if $\alg{A}$ is an algebra and $F$ is
the set of fundamental operations of~$\alg{A}$, then the
centraliser~$\cent{F}$ of~$F$ is the set
$\bigcup_{n<\omega} \Hom\apply{\alg{A}^n,\alg{A}}$. From a categorical
perspective, this is a very natural construction that makes sense in every
category~$\mathscr{C}$ with arbitrary finite powers. If~$A$ is an object in
such a category~$\mathscr{C}$ we call
$\bigcup_{n<\omega} \Hom_{\mathscr{C}}\apply{A^n,A}$ the clone over the
object~$A$. With this understanding centraliser clones are simply the
clones over algebras in the category of algebras of a certain type.
If we change the signature of the structures to allow relation symbols
(that is, we change the category to relational structures of a certain
signature), we obtain clones over some relational structure~$\mathbb{A}$
with set of fundamental relations~$Q$: $\bigcup_{n<\omega}
\Hom\apply{\mathbb{A}^n,\mathbb{A}}$. This clone is called the clone
$\Pol{Q}$ of polymorphisms of~$Q$ (or just the polymorphism clone
of the structure~$\mathbb{A}$), and it is well known by results of
Bodnar\v{c}uk, Kalu\v{z}nin, Kotov,
Romov~\cite{BodnarcukKaluzninKotovRomovGaloisTheoryForPostAlgebras} and
Geiger~\cite{GeigerClosedSystemsOfFunctionsAndPredicates}
on the classical $\PolOp$-$\InvOp$ Galois correspondence
that every clone on a finite carrier set~$A$ arises as a polymorphism
clone of some relational structure~$\mathbb{A}$.
\par

As every algebraic structure~$\alg{A}$ can also be understood as a
relational one (by taking the graphs of the fundamental operations as
the fundamental relations), it is clear that the centraliser clones on
a given set~$A$ form a subcollection of the polymorphism clones on that
set. This fact is very closely related to restricting the
$\PolOp\text{-}\InvOp$ Galois correspondence on the relational side in
such a way that the only relations taken into consideration are those
which are graphs of a function.
This restriction of the preservation relation (underlying
$\PolOp\text{-}\InvOp$) between functions and relations to functions
and function graphs leads to the notion of commutation of functions,
which is exactly the homomorphism property between finite powers of
algebras that
was used above to introduce the concept of centraliser clone. As the
Galois correspondence is restricted on one side only (the relational
one), there is a connection between the associated Galois closures: the
$\PolOp\text{-}\InvOp$ closure~$\Pol{\Inv{F}}$ of a set of
operations~$F$ (which for finite~$A$ agrees with the generated
clone~$\genClone{F}$) is weaker than the bicentrical
closure~$\bicent{F}$, that is, the double centraliser of~$F$, or,
equivalently, all functions commuting with all those functions that commute with the
functions in~$F$. The strength of the bicentrical closure in comparison
to $\Pol{\Inv{}}$ manifests itself in the following way: while
$\Pol{\Inv{F}}$ closes~$F$ against all compositions of
\nbdd{F}functions with themselves and projections (i.e.\ one iteratively
substitutes functions and variables until nothing new appears),
$\bicent{F}$ computes all functions that are primitive positively
definable from the function graphs of~$F$ (i.e.\ one interprets all
existentially quantified finite conjunctions of predicates of the form
$f(\bfa{v}) = x$ and equality predicates $y=z$ (where $f\in F$,
$\bfa{v}$ is a tuple of variables and $x,y,z$ are variables) and among
these interpretations selects those relations that are function
graphs). Functions whose graphs are constructible via such primitive
positive formul\ae{} from~$F$ have been called \emph{parametrically
expressible} through~$F$~\cite[p.~26]{KuznecovCentralisers1979}
(in contrast to functions in the clone~$\genClone{F}$ that are
\emph{explicitly expressible} via~$F$),
and also the connection of this construction with the preservation of
function graphs and the commutation of operations has first been noted
in~\cite[p.~27 et seq.]{KuznecovCentralisers1979}. For this reason
centraliser clones have also been studied under the name
\emph{parametrically closed classes}
(see e.g.~\cite{Danilcenko1978ParametricallyClosedClasses3ValuedLogic})
or~\emph{primitive positive clones}
(e.g.~\cite{BurrisWillardFinitelyManyPPClones}).
\par

It may not seem so at first glance, but the parametrical
(primitive positive, bicentrical) closure is notably
much stronger than closure under substitution. Namely, it has the
remarkable consequence that on every finite set~$A$ there are only
finitely many centraliser
clones~\cite[Corollary~4, p.~429]{BurrisWillardFinitelyManyPPClones},
which is in sharp contrast to the situation for polymorphism clones, of
which there is a continuum whenever
$\abs{A}\geq 3$~\cite{JanovMucnik1959}.
If $F$ is a centraliser clone (i.e. $\bicent{F} = F$), then
$\bicentn[1]{F}\subs\bicentn[2]{F}\subs \dotsm \subs
\bicentn{F} \subs F$ holds for all $n<\omega$ and
$\bigcup_{n<\omega} \bicentn{F} =\bicent{F} = F$. Since there are only
finitely many centraliser clones on a given finite set there must be some
$n<\omega$ such that for arities larger than~$n$ none of the inclusions is strict any more,
that is, $\bicentn{F} = \bicentn[m]{F}$ for all $n\leq
m<\omega$. Hence, $F = \bigcup_{j\leq n}\bicentn[j]{F} =
\bicentn{F}$; so there is some arity $n$ such that $F$ is
bicentrically generated by its \nbdd{n}ary part. Take this $n_F$ to be
minimal and then take the maximum over all (finitely many) $n_F$:
\[\cdeg(k)\defeq \max\lset{n_F}{F=\bicent{F} \text{ on } A,\abs{A}=k}.\]
We shall refer to this number as the \emph{uniform centraliser degree}
for a \nbdd{k}element set, since every centraliser clone~$F$ on a
carrier set of size~$k$ satisfies $F = \bicentn[\cdeg(k)]{F}$.
\par

With the help of Post's lattice, one can show that $\cdeg(2)=3$. Burris
and Willard explain in~\cite[p.~429]{BurrisWillardFinitelyManyPPClones}
that $\cdeg(k)\leq 4+k^{k^4-k^3+k^2}$ and they claim
that `[b]y slightly different methods [one] can show that any primitive
positive clone on a \nbdd{k}element set is [bicentrically] generated by
its members of arity at most~$k^k$', which implies $\cdeg(k)\leq k^k$.
No written account of the details of this argument has appeared in the
literature so far. However, at the end of the sentence cited above
Burris and Willard conjecture that $\cdeg(k)\leq k$ for every $k\geq 3$.
Besides intuition the only support for this conjecture is a series of
works by A.\,F.\ Dani\v{l}\v{c}enko on the case
$k=3$~(\cite{Danilcenko1974ParametricallyClosedClasses3ValuedLogic,
             Danilcenko1976ParametricallyIndecomposables,
             Danilcenko1977ParametricExpressibility3ValuedLogic,
             Danilcenko1978ParametricallyClosedClasses3ValuedLogic,
             Danilcenko1979-thesis},
all of these are in Russian, \cite{Danilcenko1977ParametricExpressibility3ValuedLogic}
has been translated
in~\cite{Danilcenko1977ParametricExpressibility3ValuedLogicTranslated};
\cite{Danilcenko1979ParametricalExpressibilitykValuedLogic}
is written in English). As a side note we remark that a \nbdd{k}ary
example function, stated
in~\cite[p.~269]{Danilcenko1977ParametricExpressibility3ValuedLogicTranslated}
for a different proof, can be used to show that $\cdeg(k)\geq k$ for
$k\geq 3$; so if the Burris-Willard conjecture is true, then it
certainly is sharp.
\par

In her thesis~\cite[Section~6, p.~125 et seqq.]{Danilcenko1979-thesis} Dani\v{l}\v{c}enko gives a complete description of
all $2\,986$ centraliser clones on the three\dash{}element domain. A central
step in this process is to identify a set~$\Gamma$ of 197 parametrically
indecomposable functions~\cite[Theorem~4,
p.~103]{Danilcenko1979-thesis} such that every centraliser clone~$F$ is
the centraliser of a subset of~$\Gamma$~\cite[Theorem~5,
p.~105]{Danilcenko1979-thesis}. The maximum arity of functions
in~$\Gamma$ is three, so Dani\v{l}\v{c}enko's theorems imply that
$F = \bicentn[3]{F}$ for every centraliser clone on three
elements (cf.\
Proposition~\ref{prop:char-cdeg}\eqref{item:bicent-n},\eqref{item:cent-leq-n}), that is, $\cdeg(3)\leq3$. The
results of Theorems~4 and~5 of~\cite{Danilcenko1979-thesis}, which make
the Burris-Willard conjecture true for $k=3$, are also mentioned
in~\cite[p.~155 et
seq.]{Danilcenko1979ParametricalExpressibilitykValuedLogic}
and~\cite[Section~5,
p.~414 et seq.]{Danilcenko1977ParametricExpressibility3ValuedLogic} (\cite[Section~5,
p.~279]{Danilcenko1977ParametricExpressibility3ValuedLogicTranslated},
respectively), but no proofs are given there.
\par

Drastically cut down versions of this work have been published
in~\cite{Danilcenko1974ParametricallyClosedClasses3ValuedLogic,
         Danilcenko1976ParametricallyIndecomposables,
         Danilcenko1977ParametricExpressibility3ValuedLogic,
         Danilcenko1978ParametricallyClosedClasses3ValuedLogic,
         Danilcenko1979ParametricalExpressibilitykValuedLogic},
of which
only~\cite{Danilcenko1977ParametricExpressibility3ValuedLogicTranslated,Danilcenko1979ParametricalExpressibilitykValuedLogic}
are accessible without difficulties. Given that the whole thesis
comprises~141 pages, these excerpts are rough sketches of the
classification at best (sometimes containing mistakes, many but not all
of which have been corrected in~\cite{Danilcenko1979-thesis}), and
leading experts in the field agree that it is very hard if not
impossible to reconstruct the proof of the description of all
centralisers on three\dash{}element sets from the readily available
resources. For example, Theorem~4 of~\cite{Danilcenko1979-thesis} has
appeared as part of~\cite[Proposition~2.2,
p.~16]{Danilcenko1978ParametricallyClosedClasses3ValuedLogic} with a
proof sketch of less than two pages, while the proof
from~\cite{Danilcenko1979-thesis} goes through technical calculations
and case distinctions for several pages
(however from Propositions~2.2, 2.3 and~2.4
of~\cite{Danilcenko1978ParametricallyClosedClasses3ValuedLogic}, a proof
of Theorem~5 of~\cite{Danilcenko1979-thesis} \emph{can} be obtained).
The chances of understanding might be better using the thesis as a primary source, but
for unknown reasons Moldovan librarians seem to be rather reluctant to
grant full access to it. In the light of this discussion,
Dani\v{l}\v{c}enko's classification is a result that one may believe in,
but that should not be trusted unconditionally to build further theory
on as it remains not easily verifiable at the moment. This of course
also casts some doubts on the basis of the Burris-Willard conjecture.
\par

Another possible challenge to the conjecture (and likewise to the
correctness of Dani\v{l}\v{c}enko's list of parametrically
indecomposable functions) is presented by much
later results of Snow~\cite{SnowGeneratingPrimitivePositiveClones}. In
this article the minimum arity needed to generate the bicentraliser
clone of a finite algebra from its term operations is investigated, and,
under certain assumptions on the algebra, quite satisfactory upper
bounds for that number are produced. These sometimes match (or almost
match) the number~$k$ predicted by the Burris-Willard conjecture, and
sometimes even fall below~$k$. This is possible (and supports the
conjecture) since the bounds given by Snow do not apply to \emph{all}
algebras on a \nbdd{k}element set, but only to some specific subclass.
Hence, they are not in contradiction with the \nbdd{k}ary function
from~\cite[p.~269]{Danilcenko1977ParametricExpressibility3ValuedLogicTranslated}.
Even more interestingly, in Section~3
of~\cite{SnowGeneratingPrimitivePositiveClones}
a class of examples of algebras on $k$-element carrier sets is given,
for which Snow proves $(k-1)^2$ to be a lower bound for the minimum
arity of term functions from which the bicentraliser can be generated.
This number is larger than~$k$ whenever $k\geq 3$. Explicitly,
when $k=3$, the lower bound is equal to~$4$, which means that arity
three or less does not suffice to generate the bicentraliser clone of
that specific algebra.
\par

In more detail, Snow  defines for an algebra \m{\alg{A}} with set~$F$ of
fundamental operations the number
$\ppc(\alg{A})= \min\lset{n\in\N}{ \bicent{\genClone{F}} = \bicent{{\genClone[n]{F}}}}$.
This number certainly only depends on the clone \m{\genClone{F}} of term
operations of the algebra, hence no generality is lost in simply
considering the number
\[
\mu_F \defeq
\min\lset{n\in\N}{ \bicent{\genClone{F}} = \bicent{{\genClone[n]{F}}}}
=\min\lset{n\in\N}{ \bicent{F} =\bicentn{F}}
\]
associated with clones~$F$ on a \nbdd{k}element set~$A$.
If \m{F} happens to be a centraliser clone, the definition clearly
simplifies to
\m{\mu_F = \min\lset{n\in\N}{F = \bicentn{F}} = n_F}, which is
bounded above by \m{\cdeg(k)}.
However, if now~$F$ is the clone of term operations of the example constructed
by~Snow, then the lower bound on~$\mu_F$ from~\cite[Theorem~3.1,
p.~171]{SnowGeneratingPrimitivePositiveClones} implies the
following contradiction
\[k<(k-1)^2\leq \mu_F \stackrel{?_2}{=} n_F\leq \cdeg(k)
   \stackrel{?_1}{\leq} k.\]
This offers two conclusions: either $?_1$ does not hold, which means
that the Burris-Willard conjecture and, in particular, the
Dani\v{l}\v{c}enko classification on three\dash{}element domains fail,
or $?_2$ is false, which simply means that~$F$ is not a centraliser
clone.
If we are to believe in Dani\v{l}\v{c}enko's theorems, then (for $k=3$)
the latter is the only possible consequence. However, for the reasons
mentioned above, it would be desirable to obtain such a conclusion
independently of Dani\v{l}\v{c}enko's \oe{}uvre.
\par

Such is the aim of the present article. We are going to give a proof
that the clone~$F$ of term operations of the algebra given
in~\cite[Theorem~3.1, p.~171]{SnowGeneratingPrimitivePositiveClones} is
not bicentrically closed and hence poses no threat to the Burris-Willard
conjecture. To do this, for every $k\geq 3$  we exhibit a
\nbdd{(k-1)}ary
function in~$\bicent{F}$ that cannot be obtained by composition of the
fundamental operation(s) of Snow's algebra. In doing so we use the case
$k=3$ as a guideline, where we show, for example, that~$F$
and~$\bicent{F}$ cannot be separated by unary functions, and that the
mentioned operation is the only separating binary function.

\section{Notation and preliminaries}
Throughout we use \m{\N = \set{0,1,2,\dotsc}} to denote the set of
natural numbers, and we write \m{\Np} for~$\N\setminus\set{0}$. It will
be convenient for us to understand the elements $n\in\N$ as
\nbdd{n}element sets $n= \set{0,1,\dotsc,n-1}$ as originally suggested
by John von Neumann in its model of natural numbers as finite ordinals.
\par
One of the central concepts for this paper are functions, such as
$f\colon A\to B$ and $g\colon B\to C$, and we use a
left\dash{}to\dash{}right notation for composition. That is, \m{g\circ
f\colon A\to C} sends any $x\in A$ to $g(f(x))$. The set of all
functions from~$A$ to~$B$ is written as~$B^A$.
Moreover, if \m{f\in B^A} and \m{U\subs A} and \m{V\subs B} we denote
by \m{f\fapply{U} = \set{f(x)\mid x\in U}} the \emph{image of~$U$
under~$f$} and by \m{f^{-1}\fapply{V} = \set{x\in A \mid f(x)\in V}} the
\emph{preimage of~$V$ under~$f$}. We also use the symbol~$\im f$ to
denote the full \emph{image} $f\fapply{A}$ of~$f$. All these
notational conventions will apply in particular to tuples \m{\bfa{x}\in
A^n}, \m{n\in\N}, that we formally understand as maps
\m{\bfa{x}\colon \set{0,\dotsc,n-1}\to A}. This does, of course, not
preclude us from using a different indexing for the entries
of \m{\bfa{x}=(x_1,\dotsc,x_n)}, if that seems more handy. So, e.g.,
we have \m{\im \bfa{x} = \set{x_1,\dotsc,x_n}} and \m{f\circ
\bfa{x} = (f(x_1),\dotsc, f(x_n))\in B^n}.

Notably, we are
interested in functions of the form $f\colon A^n\to A$ that we call
\emph{\nbdd{n}ary operations} on~$A$. All such operations form the set
$A^{A^n}$, and if we let the parameter~$n$ vary in~$\Np$, then we obtain
the set \m{\Op{A}= \bigcup_{0<n<\omega} A^{A^n}} of all \emph{finitary
(non\dash{}nullary) operations} over~$A$. If \m{F\subs\Op{A}} is any set of finitary
operations, we denote by \m{\Fn{F}\defeq A^{A^n}\cap F} its
\emph{\nbdd{n}ary part}. In particular, $\Op[n]{A} = A^{A^n}$.
Some specific \nbdd{n}ary operations will be needed: for $a\in A$ we
denote the constant \nbdd{n}ary function with value~$a$ by
$\cna[n]{a}\colon A^n\to A$. Moreover, if \m{n\in\N} and
\m{1\leq i\leq n} we call \m{\eni[n]{i}\colon A^n\to A}, given by
\m{\eni[n]{i}(x_1,\dotsc,x_n)\defeq x_i} for all
\m{(x_1,\dotsc,x_n)\in A^n}, the \nbdd{i}th \nbdd{n}variable
\emph{projection} on~$A$. Collecting all projections on~$A$ in one set,
we obtain
\m{\J{A} = \lset{\eni[n]{i}}{1\leq i\leq n, n\in\N}}.
\par
We call a set $F\subs\Op{A}$ a \emph{(concrete) clone} on~$A$ if
$\J{A}\subs F$ and if~$F$ is closed under composition, i.e.,
whenever \m{m,n\in\N} and \m{f\in\Fn{F}},
\m{g_1,\dotsc,g_n\in\Fn[m]{F}}, then
also the composition \m{f\circ(g_1,\dotsc,g_n)}, given by
\m{(f\circ(g_1,\dotsc,g_n))(\bfa{x}) \defeq
f(g_1(\bfa{x}),\dotsc,g_n(\bfa{x}))} for any \m{\bfa{x}\in A^m},
belongs to the set~$F$. All sets of operations that were named `clone'
in the introduction are indeed clones in this sense (except for the
fact that they were allowed to contain nullary operations, which we
want to exclude to avoid unnecessary technicalities). Clones are closed
under intersections, and hence for any set $G\subs\Op{A}$ there is a
least clone~$F$ under inclusion with the property \m{G\subs F}. This
clone~$F$ is called the \emph{clone generated by~$G$} and is denoted
as~$\genClone{G}$. It is computed by adding all projections to~$G$ and
then closing under composition, that is, by forming all term operations
(of any positive arity) over the algebra~$\algwops{A}{G}$.
\par
A function \m{f\in\Op[n]{A}} \emph{preserves} a relation
\m{\rho\subs A^m} (with \m{m,n\in\N}) if for every
\m{\bfa{r}=(r_1,\dotsc,r_n)\in\rho^n} the
tuple \m{f\circ\bfa{r}\defeq (f(r_1(i),\dotsc,r_n(i)))_{1\leq i\leq m}}
belongs to~$\rho$. For a set~$Q$ of finitary relations, the
set \m{\Pol{Q}} of polymorphisms of~$Q$ consists of all functions
preserving all relations belonging to~$Q$. Every polymorphism set is a
clone. Dually, for a set~$F\subs\Op{A}$, the set \m{\Inv{F}} contains
all \emph{invariant} relations of~$F$, that is, all relations being
preserved by all functions in~$F$.
\par
For the convenience of the reader we now give a perhaps more accessible
characterisation of the (non\dash{}nullary part of the) centraliser \m{\cent{F}} of some set of
operations~$F\subs\Op{A}$, which was already defined at the beginning
of the introduction. A function \m{g\colon A^m\to A} belongs to the
centraliser~$\cent{F}$ (\emph{commutes} with all functions from~$F$) if for every function \m{f\in F} the following
holds (where~$n$ is the arity of~$f$): for every matrix \m{X\in
A^{m\times n}} applying~$g$ to the \nbdd{m}tuple obtained from
applying~$f$ to the rows of~$X$ gives the same result as evaluating~$f$
on the \nbdd{n}tuple obtained from applying~$g$ to the columns of the
matrix. In symbols:
\m{g((f((x_{ij})_{j\in n}))_{i\in m})
 = f((g((x_{ij})_{i\in m}))_{j\in n})}
has to hold for all
\m{(x_{ij})_{(i,j)\in m\times n}\in A^{m\times  n}}
(and all $f\in F$).
A brief moment of reflection shows that this condition is the same as
saying that \m{g\colon \algwops{A}{F}^m\to \algwops{A}{F}} is a
homomorphism. A yet different way of saying this is that~$g$ is a
polymorphism of~$\mathbb{A}=\algwops{A}{\graph{F}}$, that is,
$g\in\Pol{\graph{F}}$ preserves all graphs
\m{\graph{f}=\lset{(\bfa{x},f(\bfa{x}))}{\bfa{x}\in A^n}\subs A^{n+1}}
of all functions \m{f\in F} of any arity \m{n\in \N}. From this, it is
again clear that \m{\cent{F}} always must be a clone. On the other
hand, it is obvious from the matrix formulation that centralisation is
a symmetric condition:
for all \m{F,G\subs\Op{A}} we have \m{G\subs \cent{F}} if and only if
\m{F\subs \cent{G}}. Hence, we see that
\begin{align*}
\cent{F} &= \lset{g\in\Op{A}}{g\in\cent{F}}
=\lset{g\in\Op{A}}{F\subs\cent{\set{g}}}\\
&=\lset{g\in\Op{A}}{\genClone{F}\subs\cent{\set{g}}}
=\lset{g\in\Op{A}}{g\in\cent{\genClone{F}}} = \cent{\genClone{F}}
\end{align*}
for every \m{F\subs\Op{A}}, so the centraliser of a whole clone is not
smaller than the centraliser of its generators. Since the clone
constructed in Snow's paper is given in terms of a single generator
function, we can thus study its centraliser as the set of all operations
commuting with this one generating function.
\par

In the introduction the uniform centraliser degree was defined as the
least arity~$n$ such that every centraliser clone~$F$ on a given finite
set can be bicentrically generated as \m{F=\bicentn{F}}. The
following result shows that the search for this number is likewise a
search for an arity~$n$ such that every centraliser clone is a
centraliser of a set of functions of arity at most~$n$.
\begin{proposition}\label{prop:char-cdeg}
For any carrier set~$A$ and an integer~$n\in\N$ the following facts are
equivalent:
\begin{enumerate}[(a)]
\item\label{item:bicent-n}
      For every centraliser clone~$F$ we have \m{F=\bicentn{F}}.
\item\label{item:n-cent}
      For every centraliser clone~$F$ we have
      \m{\centn{F}=\cent{F}}.
\item\label{item:n-cent-n-cent}
      For every centraliser clone~$F$ we have
      $F^{(n)*(n)*}=F$.
\item\label{item:cent-leq-n}
      For every centraliser clone~$F$ there is some
      \m{G\subs \bigcup_{\ell\leq n} \Op[\ell]{A}} such that
      \m{F=\cent{G}}.
\item\label{item:cent-n}
      For every centraliser clone~$F$ there is some
      \m{G\subs \Op[n]{A}} such that
      \m{F=\cent{G}}.
\item\label{item:cent-n-cent-bicent}
      For every set~$F\subs\Op{A}$ we have
      $F^{*(n)*}=\bicent{F}$.
\item\label{item:cent-n-cent}
      For every centraliser clone~$F$ we have
      $F^{*(n)*}=F$.
\end{enumerate}
\end{proposition}
\begin{proof}
If~\eqref{item:bicent-n} holds and~$F$ is a centraliser clone, then
\m{\cent{F}=F^{(n)***}=\centn{F}}, so~\eqref{item:n-cent} is
true. If~\eqref{item:n-cent} holds, then
\m{F=\bicent{F}=\bicentn{F}} for any centraliser clone~$F$, so
\m{\eqref{item:bicent-n}\Leftrightarrow\eqref{item:n-cent}}.
\par
Suppose now that~\eqref{item:bicent-n}, and thus~\eqref{item:n-cent},
hold. Letting \m{G\defeq\centn{F}} for a centraliser clone~$F$,
we have \m{F=\bicentn{F}=\cent{G}} from~\eqref{item:bicent-n}.
Applying now~\eqref{item:n-cent} to the centraliser~$G$ gives
\m{F=\cent{G}=\centn{G}=F^{(n)*(n)*}},
so~\eqref{item:bicent-n} implies~\eqref{item:n-cent-n-cent}.
\par
From~\eqref{item:n-cent-n-cent} we get~\eqref{item:cent-n} by letting
\m{G=F^{(n)*(n)}\subs\Op[n]{A}}, and~\eqref{item:cent-n}
directly gives~\eqref{item:cent-leq-n}.
\par
Now, suppose that~\eqref{item:cent-leq-n} holds for~$F$ with
functions~$G$ of arity at most~$n$.
Since we have excluded nullary operations, this implies that
\m{G\subs \genClone{\genClone[n]{G}}}, so we obtain
\m{\cent{G}\sups\cent{\genClone{\genClone[n]{G}}} =
\cent{{\genClone[n]{G}}}\sups\cent{\genClone{G}}=\cent{G}},
which means that \m{F=\cent{G} = \cent{H}} where
\m{H\defeq\genClone[n]{G}\subs\Op[n]{A}}. Thus
\m{\eqref{item:cent-n}\Leftrightarrow\eqref{item:cent-leq-n}}.
\par
From~\eqref{item:cent-n}, for every \m{F\subs\Op{A}}, we can express
the bicentraliser \m{\bicent{F}=\cent{G}} with some \m{G\subs\Op[n]{A}}.
Clearly, \m{G\subs\bicent{G}=\cent{F}}, so
\m{G\subs\ncent{F}\subs\cent{F}}.
Therefore, we obtain
\m{\bicent{F}=\cent{G}\sups F^{*(n)*}\sups\bicent{F}},
i.e.~\eqref{item:cent-n-cent-bicent}. The latter
entails~\eqref{item:cent-n-cent} as a special case, for every
centraliser clone~$F$ satisfies \m{\bicent{F}=F}. Moreover,
\eqref{item:cent-n-cent} directly gives~\eqref{item:cent-n} by letting
\m{G\defeq \ncent{F}\subs\Op[n]{A}}.
\par
It remains to show that~\eqref{item:cent-n-cent}
implies~\eqref{item:bicent-n}. Namely, for a centraliser clone~$F$,
applying~\eqref{item:cent-n-cent} to \m{G=\cent{F}}, we get
\m{G=G^{*(n)*}=F^{**(n)*} = \centn{F}}, so
\m{\bicentn{F}=\cent{G}=F}.
\end{proof}

\begin{remark}\label{rem:equivalences-for-one-F}
A closer inspection of the proof of Proposition~\ref{prop:char-cdeg}
reveals that for an individual centraliser clone~$F$ the conditions in
statements~\eqref{item:bicent-n} and~\eqref{item:n-cent} are equivalent
without the universal quantifier. The same holds for the
facts~\eqref{item:cent-leq-n}, \eqref{item:cent-n}
and~\eqref{item:cent-n-cent}.
\end{remark}

Let us now assume that~$F$ denotes the clone constructed by Snow
in~\cite{SnowGeneratingPrimitivePositiveClones}. It is our aim to show
that there is a separating function \m{f\in\bicent{F}\setminus F}.
Since the clone~$F$ is given
in~\cite{SnowGeneratingPrimitivePositiveClones} as
\m{F=\genClone{\set{T}}} by means of a generating function~$T$, once
we have selected an \nbdd{n}ary candidate function~$f$, it is not too
hard to show that \m{f\notin F}. One simply has to describe the
\nbdd{n}ary term operations of~$T$ and to show that~$f$ is not among
them. The harder part is to choose a suitable function
\m{f\in\bicent{F}}: by the definition of the bicentraliser one first
has to understand the whole set~$\cent{F}$ in order to
calculate~$\bicent{F}$. As~$\cent{F}$ contains functions of all arities
this task may require infinitely many steps. Admittedly, there is an
upper bound on the arities that have to be considered, but this bound
is connected to~$\cdeg(\abs{A})$ (see
\m{\eqref{item:bicent-n}\Leftrightarrow\eqref{item:cent-n-cent-bicent}}
in Proposition~\ref{prop:char-cdeg}) and hence under current knowledge
the number of steps is at least exponentially big.
\par
As a way out of this dilemma, we can however consider upper
approximations of~$\bicent{F}$. Namely, if we cut down the centraliser
at some arity~$\ell$, then \m{F^{*(\ell)*}\sups\bicent{F}}.
The smaller~$\ell$ the coarser these approximations are,
but also the easier it becomes to describe \m{\ncent[\ell]{F}}.
In the subsequent section we shall employ a strategy, where we
always start with the least interesting arity~$\ell=1$; it turns
out that this already produces good results by ruling out many
functions that cannot belong to~$\bicent{F}$.
\par

To obtain more information about \m{\ncent[\ell]{F}} for
some fixed~$\ell$, it will be important to derive as many necessary
conditions as possible to help to narrow down the possible candidate
functions in the centraliser. This is done by observing that any
\m{g\in \cent{F}=\Pol{\graph{F}}} belongs to
\m{\Pol{\Inv{\Pol{\graph{F}}}}} and thus has to
preserve all relations in the relational clone
\m{\Inv{\Pol{\graph{F}}}} generated by the graphs of the functions
from~$F$. This set contains all relations that can be defined via
primitive positive formul\ae{} from \m{\graph{F}=\lset{\graph{f}}{f\in
F}}, and among these there are a few notorious candidates: the image,
the set of fixed points and the kernel of any function \m{f\in\Fn{F}}:
\begin{align*}
\im(f)&= \lset{z\in A}{\exists x_1,\dotsc,x_n\in A\colon z =
f(x_1,\dotsc,x_n)},\\
\fix(f) &= \lset{z\in A}{f(z,\dotsc,z)=z},\\
\ker(f) &= \lset{(x_1,\dotsc,x_{2n})\in
A^{2n}}{\exists z\in A\colon f(x_1,\dotsc,x_n)=z
=f(x_{n+1},\dotsc,x_{2n})}.
\end{align*}
\par
To make this more concrete, we now give the generating
function~\m{T\in\Op[n^2]{A}} for the clone \m{F=\genClone{\set{T}}}
where \m{\mu_F\geq n^2} on \m{A=\set{0,\dotsc,n}}, \m{n\geq 2}
(see p.~172 of~\cite{SnowGeneratingPrimitivePositiveClones}):
\m{T\apply{x_{11},\dotsc,x_{1n},x_{21},\dotsc,x_{2n},\dotsc,x_{n1},\dotsc,x_{nn}}=1}
if \m{x_{ij}=i} for all \m{i,j\in\set{1,\dotsc,n}} or
\m{x_{ij}=j} for all \m{i,j\in\set{1,\dotsc,n}}, and
it is zero for all other arguments.
Hence, \m{\im(T) = \set{0,1}}, \m{\fix(T) = \set{0}} and
\m{\ker(T)} identifies
\m{(1,2,\dotsc,n,\dotsc,1,2,\dotsc,n)} with
\m{(1,1,\dotsc,1,\dotsc,n,n,\dotsc,n)} in one block, and all other
\nbdd{n^2}tuples in a second block.
\par
Eventually, after we have found a suitable candidate function
\m{f\notin F=\genClone{\set{T}}}, upper approximations
\m{F^{*(\ell)*}} will not any more be enough to prove
that \m{f\in\bicent{F}} (unless we use an exponentially high value
for~$\ell$, cf.\
Proposition~\ref{prop:char-cdeg}\eqref{item:bicent-n},\eqref{item:cent-n-cent-bicent}). Instead, we can apply a Galois theoretic trick. Namely,
\m{f\in\bicent{F}} if and only if
\m{\cent{F}\subs\cent{\set{f}}=\Pol{\set{\graph f}}}, which is
equivalent to \m{\graph f\in\Inv{\cent{F}} = \Inv{\Pol{\graph{F}}}}.
As the carrier set is finite, this means that the graph of~$f$ must
belong to the relational clone generated from the graphs of functions
in~$F$, i.e., that it is primitive positively definable from those
graphs. Finding a primitive positive formula, which does the job,
requires some creativity, and we will try our best to give some
intuition how it can be found in the case where $\abs{A}=3$. For the
general case $\abs{A}=k\geq 3$ we shall only state the generalisation
of the respective formula and verify that it suffices to define the
graph of a \nbdd{(k-1)}ary function that does not belong to~$F$.

\section{Separating a clone from its bicentraliser}
For the remainder of the paper we let \m{A=\set{0,1,\dotsc,k-1}} where
\m{k\geq 3}, and we consider the clone \m{F=\genClone{\set{T}}}
constructed by Snow
in~\cite[Section~3]{SnowGeneratingPrimitivePositiveClones}. For the
definition of the \nbdd{(k-1)^2}ary generating function~$T$, see the
end of the preceding section.
\par
It is our task to identify some arity~$n$ and some \nbdd{n}ary
operation \m{f\in\Op[n]{A}} such that
\m{f\in\bicent{F} = \bicent{\set{T}}}, but
\m{f\notin F=\genClone{\set{T}}}. In order to avoid a combinatorial
explosion of the structure of the involved clones, it is of course
desirable to keep the arity~$n$ as low as possible. Hence, we shall
start with a description of~\m{\genClone[n]{\set{T}}} for \m{n<k-1}.
Then, using the method of upper approximations, we shall show that it
is impossible to find a separating \m{f\in\bicent{\set{T}}} of such a
low arity. So the next step will be to consider $n=k-1$. Here, we will
first study the case $k=3$, where we can show that there is a unique
function of arity \m{n=k-1=2}, for which we can prove
\m{f\in\bicent{\set{T}}}, but \m{f\notin \genClone[2]{\set{T}}}.
Subsequently, we shall demonstrate that the construction of this
particular~$f$ (and the proof of \m{f\in\bicent{\set{T}}}) can be
generalised to any $k\geq 3$.

\begin{lemma}\label{lem:Tclone-small}
For any $k\geq 3$ we have
\m{\genClone[n]{\set{T}} = \J[n]{A} \cup\set{\cna[n]{0}}} for all $1\leq n<k-1$
where $\cna[n]{0}$ denotes the \nbdd{n}ary constant zero function.
\end{lemma}
\begin{proof}
We have $\cna[n]{0}= \composition{T}{\eni{1},\dotsc,\eni{1}}$ since~$T$
maps every constant tuple to~$0$. Thus the mentioned functions belong to
the \nbdd{n}ary part of~$\genClone{\set{T}}$.
Moreover, the given set is a subalgebra of \m{\algwops{A}{T}^{A^n}}: namely
every composition of~$T$ with functions at least one of which is
$\cna[n]{0}$ is $\cna[n]{0}$. This is so since~$T$ maps every tuple
containing a zero entry to zero. Furthermore, every composition of~$T$
involving only (some of) the~$n$ projections is also~$\cna[n]{0}$
as~$T$ maps every tuple with at most~$n<k-1$ distinct entries to zero.
\end{proof}

To describe \m{\nbicent{\set{T}}} for \m{0<n<k-1}, we shall study
lower approximations of~\m{\cent{\set{T}}}. We begin by cutting the
arity at the level \m{\ell=1}.

\begin{lemma}\label{lem:Tstar1}
For~$A=\set{0,\dotsc,k-1}$ of size $k\geq 3$ we have\footnote{%
For $k=3$ the correctness of this lemma can be checked with the
Z3-solver~\cite{deMouraBjoernerZ3EfficientSMTsolver,Z3} using
the ancillary file \texttt{unaryfunccommutingT.z3}.
It can also be seen from the file \texttt{Tcent1.txt} produced by the
function \texttt{findallunaries()} from the ancillary file
\texttt{commutationTs.cpp}.}
\[\ncent[1]{\set{T}}
               = \set{\id_A} \cup \lset{f\in\Op[1]{A}}{f(0)=f(1)=0}.\]
\end{lemma}
\begin{proof}
Let us fix \m{f\in\Op[1]{A}}, commuting with \m{T}.
Since \m{f\in \Pol{\set{\fix(T)}}}, we have \m{f(0) = 0}. Moreover, since $f$
preserves the image of~$T$, we must have \m{f(1)\in\set{0,1}}. If
$f(1)=0$, we are done. Otherwise, if $f(1)=1$, we shall show that
$f=\id_A$. Namely, since $f$ and~$T$ commute, we have
\begin{align*}
1= f(1) &= f(T(1,\dotsc,1,2,\dotsc,2,\dotsc,k-1,\dotsc,k-1))\\
   &=T(f(1),\dotsc,f(1),f(2),\dotsc,f(2),\dotsc,f(k-1),\dotsc,f(k-1)),
\end{align*}
which implies that
\begin{align*}
(f(1),\dotsc,f(1),&f(2),\dotsc,f(2),\dotsc,f(k-1),\dotsc,f(k-1))\\
&\in
T^{-1}[\set{1}]\setminus\set{(1,2,\dotsc,k-1,1,2,\dotsc,k-1,\dotsc,1,2,\dotsc,k-1)}\\
&=\set{(1,\dotsc,1,2,\dotsc,2,\dotsc,k-1,\dotsc,k-1)},
\end{align*}
whence clearly $f(x) = x$ for all $0<x<k$, i.e.\ \m{f=\id_A}.
\par
Conversely, we prove that every $f\in\Op[1]{A}$ with $f(0)=f(1)=0$
commutes with~$T$. Assume, for a contradiction, that for some
\m{\bfa{x}\in A^{(k-1)^2}} we had $T(f\circ\bfa{x})=1$; then
\m{\set{1,\dotsc,k-1} = \im f\circ\bfa{x} \subs \im f}, so~$f$ would be
surjective, and, by finiteness of~$A$, bijective. This would contradict
$f(0)=f(1)=0$, so for every \m{\bfa{x}\in A^{(k-1)^2}} we have
$T(f\circ \bfa{x}) = 0 = f(0) = f(1) = f(T(\bfa{x}))$,
since \m{T(\bfa{x})\in \set{0,1}}. Thus~$f\in\cent{\set{T}}$.
\end{proof}

\begin{corollary}\label{cor:Tstar1}
For $A=\set{0,\dotsc,k-1}$ of cardinality~$k\geq 3$, we have the inclusion
\[\ncent[1]{\set{T}}
   \sups \lset{u_{j,a}}{a\in A\land j\in A\setminus\set{0,1}},\]
where \m{u_{j,a}} is given by the rule
\[u_{j,a}(x) = \begin{cases} a&\text{if }x=j,\\
                             0&\text{otherwise.}
               \end{cases}\]
\end{corollary}

The binary part of the centraliser already becomes rather obscure
in the general case. So we only give a description for
the case $k=3$ (which can certainly also be verified by a brute-force
enumeration using a computer).

\begin{lemma}\label{lem:Tstar2}
For \m{A=\set{0,1,2}} the set \m{\ncent[2]{\set{T}}} contains the following 65 functions
\begin{align*}
\ncent[2]{\set{T}}= \set{\eni[2]{1},\eni[2]{2}}
&{}\disjointunion \lset{z_a}{a\in \set{0,1,2}}\\
&{}\disjointunion \bigcup_{c\in\set{1,2}}
\lset{f_{a,\bfa{x}}}{a\in\set{0,c}\land \bfa{x}\in\set{0,c}^4\setminus\set{(0,0,0,0)}}
\end{align*}
given by the following tables\footnote{%
The correctness of this lemma (and its proof) can be checked with the
Z3-solver~\cite{deMouraBjoernerZ3EfficientSMTsolver,Z3} using the
ancillary file \texttt{binaryfunccommutingT.z3}. The completeness of
the list of 65 operations can also be verified with the function
\texttt{findallbinaries()} from the ancillary file
\texttt{commutationTs.cpp}, resulting in the file \texttt{Tcent2.txt}.}:
\begin{align*}
&\begin{array}{r|*{3}{c}}
z_a(x\backslash y)& 0&1&2\\\hline
0& 0& 0& 0\\
1& 0& 0& 0\\
2& 0& 0& a
\end{array}&
&\begin{array}{r|*{3}{c}}
f_{a,(b,c,d,e)}(x\backslash y)& 0&1&2\\\hline
0& 0& 0& b\\
1& 0& 0& c\\
2& d& e& a
\end{array}
\end{align*}
\end{lemma}
\begin{proof}
Using a case distinction, one can verify that every function
\m{f\in \ncent[2]{\set{T}}} must be among the ones mentioned in
the lemma.
\begin{enumerate}[1.]
\item Assume \m{f(1,1) = 1}. We can show that \m{f(2,2)=2} and
      \m{\set{f(1,2),f(2,1)} = \set{1,2}}. Namely, \m{f\in\cent{\set{T}}}
      implies
      \begin{multline*}
      1=f(1,1) = f(T(1,1,2,2),T(1,2,1,2))\\
       =T(f(1,1),f(1,2),f(2,1),f(2,2)) = T(1,f(1,2),f(2,1),f(2,2)),
      \end{multline*}
      which is only possible if \m{f(2,2)=2} and
      \m{(f(1,2),f(2,1))\in\set{(1,2),(2,1)}}.
\begin{enumerate}[{1.}1.]
\item Assume \m{f(1,2) =1} and \m{f(2,1) =2}. It follows that
      \m{f=\eni[2]{1}}.
      In fact, our assumption \m{f\in\cent{\set{T}}} implies
      \begin{multline*}
      f(1,0) = f(T(1,1,2,2),T(1,1,1,2))\\
       =T(f(1,1),f(1,1),f(2,1),f(2,2)) = T(1,1,2,2)=1,
      \end{multline*}
      moreover
      \begin{multline*}
      f(0,1) = f(T(1,1,1,2),T(1,1,2,2))\\
       =T(f(1,1),f(1,1),f(1,2),f(2,2)) = T(1,1,1,2)=0,
      \end{multline*}
      and
      \begin{multline*}
      1 = f(1,0) = f(T(1,2,1,2),T(1,0,1,2))\\
       =T(f(1,1),f(2,0),f(1,1),f(2,2)) = T(1,f(2,0),1,2),
      \end{multline*}
      which is only possible if \m{f(2,0) = 2}.
      Finally, we have
      \begin{multline*}
      0 = f(0,1) = f(T(1,0,1,2),T(1,2,1,2))\\
       =T(f(1,1),f(0,2),f(1,1),f(2,2)) = T(1,f(0,2),1,2),
      \end{multline*}
      which means \m{f(0,2)\neq 2}, and
      \begin{multline*}
      0 = f(0,0) = f(T(0,2,1,2),T(2,2,1,2))\\
       =T(f(0,2),f(2,2),f(1,1),f(2,2)) = T(f(0,2),2,1,2),
      \end{multline*}
      which gives \m{f(0,2)\neq 1}.
      Thus \m{f(0,2)\in A\setminus\set{1,2}=\set{0}}.
\item Assume \m{f(1,2) =2} and \m{f(2,1) =1}. It follows that
      \m{f=\eni[2]{2}} by a dual argument.
      In fact, \m{f\in\cent{\set{T}}} implies
      \begin{multline*}
      f(1,0) = f(T(1,1,2,2),T(1,1,1,2))\\
       =T(f(1,1),f(1,1),f(2,1),f(2,2)) = T(1,1,1,2)=0,
      \end{multline*}
      moreover
      \begin{multline*}
      f(0,1) = f(T(1,1,1,2),T(1,1,2,2))\\
       =T(f(1,1),f(1,1),f(1,2),f(2,2)) = T(1,1,2,2)=1,
      \end{multline*}
      and
      \begin{multline*}
      1 = f(0,1) = f(T(1,0,1,2),T(1,2,1,2))\\
       =T(f(1,1),f(0,2),f(1,1),f(2,2)) = T(1,f(0,2),1,2),
      \end{multline*}
      which is only possible if \m{f(0,2) = 2}.
      Finally, we have
      \begin{multline*}
      0 = f(1,0) = f(T(1,2,1,2),T(1,0,1,2))\\
       =T(f(1,1),f(2,0),f(1,1),f(2,2)) = T(1,f(2,0),1,2),
      \end{multline*}
      which means \m{f(2,0)\neq 2}, and
      \begin{multline*}
      0 = f(0,0) = f(T(2,2,1,2),T(0,2,1,2))\\
       =T(f(2,0),f(2,2),f(1,1),f(2,2)) = T(f(2,0),2,1,2),
      \end{multline*}
      which gives \m{f(2,0)\neq 1}.
      Thus \m{f(2,0)\in A\setminus\set{1,2}=\set{0}}.
\end{enumerate}
\item Now assume that \m{f(1,1)\neq 1}. Since \m{f\in\Pol{\im(T)}}, we
      must have \m{f(1,1)=0}.
      We can show that \m{f(0,1) = 0= f(1,0)}.
      In point of fact, we have
      \begin{multline*}
      f(0,1) = f(T(1,1,0,0),T(1,1,2,2))\\
       =T(f(1,1),f(1,1),f(0,2),f(0,2)) = T(0,0,f(0,2),f(0,2)) = 0,
      \end{multline*}
      and for \m{f(1,0)=0} we argue by swapping the arguments of~\m{f}.
      \par
      Moreover, if \m{\set{1,2}\subs\set{f(0,2),f(2,0),f(1,2),f(2,1)}},
      then \m{f\notin\cent{\set{T}}}.
      Indeed, if there are \m{x,y\in \set{0,1}} such that
      \begin{enumerate}[(a)]
      \item \m{f(2,x) = 1}, \m{f(2,y)=2}, then
            \begin{multline*}
            f(T(2,2,2,2),T(x,x,y,y))= f(0,z)=0 \\
            \neq 1 =T(1,1,2,2) =T(f(2,x),f(2,x),f(2,y),f(2,y)),
            \end{multline*}
             where \m{z\in \set{0,1}}.
      \item \m{f(x,2) = 1}, \m{f(y,2)=2}, then we argue with swapped
            arguments for \m{f}.
      \item \m{f(2,x) = 1}, \m{f(y,2)=2}, then
            \begin{multline*}
            f(T(2,2,y,y),T(x,x,2,2)) = f(0,z) = 0\\
             \neq 1 = T(1,1,2,2) = T(f(2,x),f(2,x),f(y,2),f(y,2)),
            \end{multline*}
             where \m{z\in \set{0,1}}.
      \item \m{f(x,2) = 1}, \m{f(2,y)=2}, then we argue with swapped
            arguments for~\m{f}.
      \end{enumerate}
      Hence, we know that
      \m{\set{1,2}\not\subs\set{f(0,2),f(2,0),f(1,2),f(2,1)}} for
      \m{f\in\cent{\set{T}}}.
\begin{enumerate}[{2.}1.]
\item Suppose that \m{f(2,2)=0}. There is nothing left to prove: we
      already have
      \m{f\in\set{z_0}\cup\bigcup_{c\in\set{1,2}}
                       \lset{f_{0,\bfa{x}}}{\bfa{x}\in\set{0,c}^4\setminus\set{\bfa{0}}}}.
\item Suppose that \m{f(2,2)=c\in\set{1,2}} and let \m{d} be such that
      \m{\set{c,d}=\set{1,2}}.
      We prove that \m{d\notin \set{f(0,2),f(2,0),f(1,2),f(2,1)}},
      as otherwise \m{f\notin\cent{\set{T}}}. This demonstrates that
      \m{\set{f(0,2),f(2,0),f(1,2),f(2,1)}\subs\set{0,c}}, so we have
      \m{f\in\set{z_c}\cup\bigcup_{c\in\set{1,2}}
                       \lset{f_{c,\bfa{x}}}{\bfa{x}\in\set{0,c}^4\setminus\set{\bfa{0}}}}.
      \par
      For a contradiction suppose that there is some argument \m{x\in \set{0,1}}
      such that \m{f(x,2)=d}. Then for some \m{z\in\set{0,1}} we have
      \begin{multline*}
      f(T(x,x,2,2),T(2,2,2,2)) = f(z,0) = 0\\
      \neq 1 = T(d,d,c,c) =T(f(x,2),f(x,2),f(2,2),f(2,2)),
      \end{multline*}
      when \m{(c,d)=(2,1)}, and
      \begin{multline*}
      f(T(2,2,x,x),T(2,2,2,2)) = f(z,0) = 0\\
      \neq 1 = T(c,c,d,d) =T(f(2,2),f(2,2),f(x,2),f(x,2)),
      \end{multline*}
      when \m{(c,d)=(1,2)}.
      In the case where \m{f(2,x)=d} for some \m{x\in\set{0,1}} we argue
      similarly, by swapping the arguments of~\m{f}.
\end{enumerate}
\end{enumerate}
\par
For the converse inclusion, we have to check that all mentioned
functions commute with~$T$. So let $g=z_a$ for some $a\in A$ or
$g=f_{a,(b,c,d,e)}$ and consider
$x_1,\dotsc,x_4,y_1,\dotsc,y_4\in A$ to verify that~$g$ commutes
with~$T$. Put $u\defeq T(x_1,\dotsc,x_4)$ and
$v\defeq T(y_1,\dotsc,y_4)$.
Since \m{(u,v)\in \im(T)^2 = \set{0,1}^2}, we have $g(u,v) = 0$. On
the other hand, the values $w_i\defeq g(x_i,y_i)$ for $1\leq i\leq 4$
belong to $\im(g)\subs \set{0,a,b,c,d,e}$. If at least one of them
equals~$0$, then $T(w_1,\dotsc,w_4)=0$ as needed.
Otherwise, all of them belong to \m{\set{a,b,c,d,e}\setminus\set{0}}.
If $g=z_a$, then they are all equal to~$a$ and we thus have
$T(w_1,\dotsc,w_4)=0$, too. In the case that $g=f_{a,(b,c,d,e)}$, we
know from the definition of~$g$ that \m{\set{a,b,c,d,e}\subs\set{0,j}}
for some \m{j\in \set{1,2}}. Thus, $w_1=\dots =w_4=j$, and again
$T(w_1,\dotsc,w_4)=0$. In any case, we have shown
\m{g\in\cent{\set{T}}}.
\end{proof}

Next, with the help of the coarse approximations from
Lemma~\ref{lem:Tstar1}, we observe that the bicentraliser of~$T$ only
contains functions that are close to being conservative and have many
congruences.
\begin{lemma}\label{lem:almost-conservative}
For $A=\set{0,\dotsc,k-1}$ of size \m{k\geq 3} we have
\begin{multline*}
\genClone{\set{T}}\subs \bicent{\set{T}} \subs \set{T}^{*(1)*}\\
{}\subs \Pol{\lset{U\subs A}{0\in U}}
\cap\Pol{\lset{\theta\in\Eq(A)}{(0,1)\in\theta}},
\end{multline*}
where $\Eq(A)$ denotes the set of all equivalence relations on~$A$.
\end{lemma}
\begin{proof}
It is clear that \m{\set{T}^{*(1)*}
\subs\Pol{\lset{\im(f)}{f\in \ncent[1]{\set{T}}}}}
since the image of a function is primitive positively definable from its
graph. If \m{0\in U\subsetneq A}, then~$U$ contains \m{t<k-1} elements
distinct from~$0$. According to the description of the functions
in~\m{\ncent[1]{\set{T}}} given in Lemma~\ref{lem:Tstar1}, there is some
\m{f\in \ncent[1]{\set{T}}} whose image is~\m{U}.
\par
Likewise we have \m{\set{T}^{*(1)*}
\subs\Pol{\lset{\ker(f)}{f\in \ncent[1]{\set{T}}}}}
since the kernel of a function is primitive positively definable from its
graph. Any partition of~$A$ having a class containing the set \m{\set{0,1}}
can again be realised as the kernel of a function \m{f\in
\ncent[1]{\set{T}}}
since the value $f(x)$ can be chosen arbitrarily for every
\m{x\in A\setminus\set{0,1}}.
\end{proof}

Based on this lemma we can show that the \nbdd{n}ary part of
the bicentraliser of~$T$ is not bigger than \m{\genClone[n]{\set{T}}}
when $n<k-1$.
\begin{lemma}\label{lem:T*1*-small}
For~$k=\abs{A}\geq 3$ we have
\m{\nbicent{\set{T}} = \set{T}^{*(1)*(n)}= \J[n]{A} \cup\set{\cna[n]{0}}} for all $1\leq n<k-1$
where $\cna[n]{0}$ denotes the \nbdd{n}ary constant zero function.
\end{lemma}
\begin{proof}
We shall prove that
\m{\set{T}^{*(1)*(n)}\subs \J[n]{A} \cup\set{\cna[n]{0}} = \genClone[n]{\set{T}}}
(cf.\ Lemma~\ref{lem:Tclone-small}). From this it will follow
that \m{\genClone[n]{\set{T}}\subs \nbicent{\set{T}} \subs
\set{T}^{*(1)*(n)}\subs\genClone[n]{\set{T}}}
since \m{\ncent[1]{\set{T}} \subs\cent{\set{T}}} is always true.
\par
Given that~\m{A} has size \m{k\geq 3}, the set \m{A\setminus\set{0,1}}
has \m{k-2\geq n} distinct values.
Let \m{f\in \set{T}^{*(1)*(n)}}.
By Lemma~\ref{lem:almost-conservative} we know that
\m{f\in\Pol{\set{0,2\dotsc,n+1}}}, so we obtain
\m{b\defeq f(2,\dotsc,n+1)\in\set{0,2,\dotsc,n+1}}.
Now for any \m{(a_1,\dotsc,a_n)\in A^{n}} we consider the unary map~\m{u} sending
\m{j \mapsto a_{j-1}} for \m{2\leq j\leq n+1} and \m{j\mapsto 0} otherwise.
Since \m{u\in\ncent[1]{\set{T}}} by Lemma~\ref{lem:Tstar1}, we have
\m{f\in \cent{\set{u}}} and thus
\[f(a_1,\dotsc,a_{n}) = f(u(2),\dotsc,u(n+1)) = u(f(2,\dotsc,n+1))=u(b).\]
If \m{b=0}, then \m{f(a_1,\dotsc,a_n)=u(b)=0}, so \m{f=\cna[n]{0}}.
If \m{b\neq 0}, then it follows that \m{2\leq b\leq n+1}.
Thus,  we have $f(a_1,\dotsc,a_{n})=u(b)=a_{b-1}$
for all \m{(a_1,\dotsc,a_{n})\in A^n}, which shows that \m{f=\eni{b-1}}.
\end{proof}

According to Lemmata~\ref{lem:Tclone-small} and~\ref{lem:T*1*-small},
it is impossible to find
\m{f\in\nbicent{\set{T}}\setminus\genClone[n]{\set{T}}} for
\m{n<k-1} where \m{k=\abs{A}}. Next, we thus turn our attention to
\m{n=k-1}, where we will first describe \m{\genClone[k-1]{\set{T}}}:
besides projections the \nbdd{(k-1)}ary part
of~\m{\genClone{\set{T}}} contains only functions being zero everywhere
with a possible exception in only one argument tuple which may be sent
to one.
After that we shall focus for a while on the case \m{k=3} to develop
the right ideas in connection with
\m{\nbicent[k-1]{\set{T}}=\nbicent[2]{\set{T}}},
which can eventually be generalised to any \m{k\geq 3}.

\begin{lemma}\label{lem:Tclone2-general}
Given a set~$A$ of cardinality~$k\geq 3$, put $n=k-1$.
We then have
$\genClone[n]{\set{T}}=\J[n]{A}\cup\set{\cna[n]{0}}\cup F$
where $F\subs\Op[n]{A}$ is the set of \nbdd{n}ary functions in
$\genClone{\set{T}}$ which map exactly one \nbdd{n}tuple to~$1$ and everything
else to~$0$.
\end{lemma}
\begin{proof}
We have $\cna[n]{0}=
\composition{T}{\eni{1},\dotsc,\eni{1}}\in\genClone[n]{\set{T}}$ as in
Lemma~\ref{lem:Tclone-small}, so the inclusion
$G\defeq \J[n]{A}\cup \set{\cna[n]{0}}\cup F\subs\genClone[n]{\set{T}}$ is
clear. For the opposite inclusion, we prove that~$G$ is a subuniverse of
\m{\algwops{A}{T}^{A^n}}. The first step is to check that any variable
identification of~$T$ with at most~$n$ variables ends up in~$F\cup\set{\cna[n]{0}}$.
\par
Let $i\colon\set{1,\dotsc,n^2}\to\set{1,\dotsc,n}$, $j\mapsto i_j$ be a
map describing an \nbdd{n}variable identification
\m{f = \composition{T}{\eni{i_1},\dotsc,\eni{i_{n^2}}}} of~$T$.
Clearly, $\im(f)\subs\im(T)=\set{0,1}$, so every tuple that is not
mapped to one by~$f$ will be sent to zero. To obtain a contradiction,
let us assume that $\abs{f^{-1}\fapply{\set{1}}}\geq 2$. So there are
tuples $\bfa{x}\neq \bfa{y}\in A^n$ such that
\m{T\apply{x_{i_1},\dotsc,x_{i_{n^2}}} = 1 =
   T\apply{y_{i_1},\dotsc,y_{i_{n^2}}}}. The preimage
$T^{-1}\fapply{\set{1}}$ contains only two tuples, and these mention~$n$
distinct elements. To obtain one of them in the form
$\apply{x_{i_1},\dotsc,x_{i_{n^2}}}$ or
$\apply{y_{i_1},\dotsc,y_{i_{n^2}}}$ one has to use at least~$n$
distinct variable indices, so the map~$i$ has to be surjective. It is
therefore impossible that the distinct tuples~$\bfa{x}$ and~$\bfa{y}$
produce the same tuple
\m{\apply{x_{i_1},\dotsc,x_{i_{n^2}}} =
   \apply{y_{i_1},\dotsc,y_{i_{n^2}}} \in T^{-1}\fapply{\set{1}}}.
This means, one of them, say~$\bfa{x}$, gives
\m{\apply{x_{i_1},\dotsc,x_{i_{n^2}}}
  =\apply{1,\dotsc,n,1,\dotsc,n,\dotsc,1,\dotsc,n}}, from which it
follows that $i_1,\dotsc,i_n$ are all distinct (so
\m{\set{i_1,\dotsc,i_n} = \set{1,\dotsc,n}}); the other one however
produces
\m{\apply{y_{i_1},\dotsc,y_{i_{n^2}}}
  =\apply{1,\dotsc,1,2,\dotsc,2,\dotsc,n,\dotsc,n}}.
This implies that $\set{y_1,\dotsc,y_n}
=\set{y_{i_1},\dotsc,y_{i_n}}=\set{1}$, so
$\apply{y_{i_1},\dotsc,y_{i_{n^2}}} = (1,\dotsc,1)$, which is a
contradiction for $n\geq 2$.
\par
To prove that~$G$ is closed under application of~$T$, we take functions
$f_1,\dotsc,f_{n^2}$ from~$G$ and show that
$f=\composition{T}{f_1,\dotsc,f_{n^2}}\in G$.
If $f_1,\dotsc,f_{n^2}\in\J{A}$, then the composition is a variable
identification of~$T$ that belongs to~$F\cup\set{\cna[n]{0}}\subs G$.
Otherwise, suppose that (for some $1\leq j\leq n^2$) $f_j$ is a
non\dash{}projection in $F\cup\set{\cna[n]{0}}$. For every $\bfa{x}\in
A^n$ with possibly one exception we have $f_j(\bfa{x})=0$. So for all
those arguments $\bfa{x}\in A^n$, the \nbdd{j}th component of
\m{\apply{f_1(\bfa{x}),\dotsc,f_{n^2}(\bfa{x})}} contains a zero, whence
this tuple is mapped to zero by~$T$. Consequently, $f(\bfa{x})=0$ for
all but possibly one $\bfa{x}\in A^n$, so
\m{f\in F\cup\set{\cna[n]{0}}\subs G}.
\end{proof}

For three\dash{}element domains we obtain a more specific result.
\begin{lemma}\label{lem:Tclone2}
For $A=\set{0,1,2}$ we have \m{\genClone[2]{\set{T}} = \set{\eni[2]{1},\eni[2]{2}, \cna[2]{0}, \delta_{(1,2)}, \delta_{(2,1)}}},
where \m{\cna[2]{0}} is the constant zero function and \m{\delta_a(x) = 1}
if \m{x=a} and \m{\delta_a(x)=0} otherwise.
\end{lemma}
\begin{proof}
It is easy to see that the listed binary functions belong to the clone, namely
\begin{align*}
\cna[2]{0}     &= \composition{T}{\eni[2]{1},\eni[2]{1},\eni[2]{1},\eni[2]{1}},\\
\delta_{(1,2)} &= \composition{T}{\eni[2]{1},\eni[2]{1},\eni[2]{2},\eni[2]{2}}
                = \composition{T}{\eni[2]{1},\eni[2]{2},\eni[2]{1},\eni[2]{2}},\\
\delta_{(2,1)} &= \composition{T}{\eni[2]{2},\eni[2]{1},\eni[2]{2},\eni[2]{1}}
                = \composition{T}{\eni[2]{2},\eni[2]{2},\eni[2]{1},\eni[2]{1}}.
\end{align*}
It is not hard to verify that the given subset is a subuniverse of
\m{\algwops{A}{T}^{A^2}}. Any \nbdd{T}composition involving only
projections except for the ones shown to yield~\m{\delta_{(1,2)}}
or~\m{\delta_{(2,1)}} produces~\m{\cna[2]{0}}. Any composition
involving~\m{\cna[2]{0}}, or just the \m{\delta_a}~functions yields again
the constant map~\m{\cna[2]{0}}. Therefore, only compositions involving
the \m{\delta_a}~functions \emph{and} projections have to be checked. If
all four of them are substituted into~\m{T} (in any order), the result
is~\m{\cna[2]{0}}. If only one projection (and possibly some
non\dash{}projections) are substituted, then in most cases, the result
is~\m{\cna[2]{0}}, and for a few substitutions it is one of the
\m{\delta_a}~functions. If both projections and only one of the
\m{\delta_a}~functions are substituted, the result is either the
substituted function~\m{\delta_a} or~\m{\cna[2]{0}}.
\end{proof}

With the aim of finding separating binary functions in
\m{\bicent{\set{T}}} for \m{\abs{A}=3}, we collect some properties of
binary operations in upper approximations of~\m{\bicent{\set{T}}}.

\begin{lemma}\label{lem:observations-T*1*}
Let \m{A=\set{0,1,2}} and
\m{g\in \set{T}^{*(1)*(2)}},
then the following implications hold:
\begin{enumerate}[(a)]
\item \m{g(1,2)=2 \implies \forall a\in A\colon g(0,a)=a}.
\item \m{g(2,1)=2 \implies \forall a\in A\colon g(a,0)=a}.
\item \m{g(1,2)\in\set{0,1} \implies \forall a\in A\colon g(0,a)=0}.
\item \m{g(2,1)\in\set{0,1} \implies \forall a\in A\colon g(a,0)=0}.
\end{enumerate}
\end{lemma}
\begin{proof}
By Corollary~\ref{cor:Tstar1} we have \m{g\in\cent{\lset{u_{2,a}}{a\in
A}}}. This implies for all \m{a\in A} that
\m{a=u_{2,a}(2)=u_{2,a}(g(1,2)) = g(u_{2,a}(1),u_{2,a}(2))=g(0,a)}
provided~\m{g(1,2)=2}. A symmetric argument works for \m{g(2,1)=2}.
Similarly, if \m{g(1,2)\in\set{0,1}}, then
\m{0=u_{2,a}(g(1,2)) = g(u_{2,a}(1),u_{2,a}(2))=g(0,a)}
for all \m{a\in A}, and symmetrically, if \m{g(2,1)\in\set{0,1}}.
\end{proof}

Not very surprisingly, \m{\ncent[1]{\set{T}}} does not encode
enough information about \m{\cent{\set{T}}} to determine functions in
\m{\bicent{\set{T}}} sufficiently well. However, using the description
of \m{\ncent[2]{\set{T}}} available for \m{\abs{A}=3} from
Lemma~\ref{lem:Tstar2}, we are able to derive a more promising result:
for \m{\abs{A}=3} there is a unique binary function in
\m{\set{T}^{*(2)*(2)}\setminus\genClone{\set{T}}}.
This function might---and although we do not know it yet at this point,
it actually will---serve to distinguish \m{\bicent{\set{T}}} and
\m{\genClone{\set{T}}}.
\begin{lemma}\label{lem:unique-bin-func}
For \m{A=\set{0,1,2}} we have\footnote{%
The correctness of this lemma can be checked with the
Z3-solver~\cite{deMouraBjoernerZ3EfficientSMTsolver,Z3} using the
ancillary file \texttt{func\_Tc2c2.z3}.}
\m{\set{T}^{*(2)*(2)}=\genClone[2]{\set{T}}\mathbin{\dot{\cup}}\set{f}}
where for all \m{x,y\in A}
\[f(x,y) = \begin{cases}
  1 &\text{if } \set{x,y} = \set{1,2},\\
  0 &\text{else.}
  \end{cases}\]
\end{lemma}
\begin{proof}
The proof is by a systematic case distinction.
Let \m{g\in\set{T}^{*(2)*(2)}}, which implies that
\m{g\in\set{T}^{*(2)*}=\cent{\genClone{\ncent[2]{\set{T}}}}
\subs \set{T}^{*(1)*}} since
\m{\ncent[1]{\set{T}}\subs\genClone{\ncent[2]{\set{T}}}}.
Hence, we can apply the implications from Lemma~\ref{lem:observations-T*1*} to~$g$.
\par
Assume \m{g(1,2) = 2}.
It follows by Lemma~\ref{lem:observations-T*1*} that
\m{g(0,a)=a} for all \m{a\in A}.
Our goal is to show that \m{g=\eni[2]{2}}.
For a contradiction, suppose that \m{g(2,1)=2}. Since
\m{g\in\cent{\set{z_1}}} by Lemma~\ref{lem:Tstar2}, we obtain
\[
1 = z_1(2,2) = z_1(g(1,2),g(2,1)) = g(z_1(1,2),z_1(2,1))
  = g(0,0),
\]
in contradiction to $g(0,0)=0$ derived above.
Hence \m{g(2,1)\in\set{0,1}}. Using again
Lemma~\ref{lem:observations-T*1*}, this implies \m{g(a,0)=0}
for all \m{a\in A}.
\par
Again, for a contradiction, we suppose that \m{g(2,1)=0}. Since
\m{g\in\cent{\set{f_{0,(1,1,1,0)}}}} by Lemma~\ref{lem:Tstar2}, we get
\begin{align*}
1&=f_{0,(1,1,1,0)}(2,0)=f_{0,(1,1,1,0)}(g(1,2),g(2,1))\\
 &= g(f_{0,(1,1,1,0)}(1,2),f_{0,(1,1,1,0)}(2,1)) = g(1,0),
\end{align*}
which contradicts \m{g(1,0)=0}.
\par
Hence \m{g(2,1)=1}. Then, since
\m{g\in\cent{\set{f_{0,(c,c,c,c)}}}} for \m{c\in\set{1,2}} by
Lemma~\ref{lem:Tstar2}, we get
\begin{align*}
c&=f_{0,(c,c,c,c)}(2,1)=f_{0,(c,c,c,c)}(g(1,2),g(2,1))\\
 &= g(f_{0,(c,c,c,c)}(1,2),f_{0,(c,c,c,c)}(2,1)) = g(c,c),
\end{align*}
which shows that \m{g=\eni[2]{2}}.
Note that a symmetric argument shows that the assumption
\m{g(2,1)=2} implies \m{g=\eni[2]{1}}.
\par
Assume \m{\set{g(1,2),g(2,1)}\subs\set{0,1}}. By
Lemma~\ref{lem:observations-T*1*} we get that
\m{g(0,a)=g(a,0)=0} for all \m{a\in A}. Clearly,
\m{h=\composition{g}{\id_A,\id_A}\in
\set{T}^{*(2)*(1)}\subs\set{T}^{*(1)*(1)}= \set{\id_A,\cna[1]{0}}},
see Lemma~\ref{lem:T*1*-small}.
For a contradiction, suppose that \m{h=\id_A},
whence \m{g(2,2)=2}. As \m{g\in\cent{\set{f_{0,(2,2,2,2)}}}} by
Lemma~\ref{lem:Tstar2}, we get
\begin{align*}
      2&=f_{0,(2,2,2,2)}(2,g(1,2))=f_{0,(2,2,2,2)}(g(2,2),g(1,2))\\
       &= g(f_{0,(2,2,2,2)}(2,1),f_{0,(2,2,2,2)}(2,2)) = g(2,0),
\end{align*}
which contradicts \m{g(2,0)=0} from before. Therefore,
\m{h=\cna[1]{0}}, which shows that
\m{g\in\set{\cna[2]{0},\delta_{(1,2)},\delta_{(2,1)},f}}.
\par
Hence, according to Lemma~\ref{lem:Tclone2}, \m{g\in\genClone[2]{\set{T}}\cup\set{f}}. For the converse inclusion,
one uses that the containment
\m{\genClone{\set{T}}\subs\bicent{\set{T}}\subs\set{T}^{*(2)*}}
is trivially true and one verifies that, indeed,
\m{f\in\set{T}^{*(2)*}}.
We postpone the latter until Lemma~\ref{lem:pp-formula}, where we shall
show more generally that even
\m{f\in\bicent{\set{T}}\subs\set{T}^{*(2)*}}.
Alternatively, one may ask a computer to check that~$f$ commutes with
all the $65$~functions given in Lemma~\ref{lem:Tstar2},
immediately giving a positive answer.\footnote{%
This can, for example, be done with the
Z3-solver~\cite{deMouraBjoernerZ3EfficientSMTsolver,Z3} using the
ancillary file \texttt{func\_Tc2c2.z3}.}
\end{proof}

So far, for the binary operation~$f$ exhibited in
Lemma~\ref{lem:unique-bin-func} we do not know whether it actually
belongs to~\m{\bicent{\set{T}}} as we have only worked with upper
approximations of this bicentraliser, not with~\m{\bicent{\set{T}}}
itself.
\begin{remark}\label{rem:f-in-cent-3-cent-T}
Without much more ingenuity but some additional computational effort,
it is possible to show that the unique binary operation~$f$ from
Lemma~\ref{lem:unique-bin-func} belongs to
\m{\set{T}^{*(3)*}}, which is even closer to \m{\bicent{\set{T}}}.
\par
To do this one needs to enumerate~\m{\ncent[3]{\set{T}}}. Since
\m{\cent{\set{T}}} is a clone, for every ternary $g\in\cent{\set{T}}$
each of its identification minors
$\composition{g}{\eni[2]{1},\eni[2]{1},\eni[2]{2}}$,
$\composition{g}{\eni[2]{1},\eni[2]{2},\eni[2]{1}}$ and
\m{\composition{g}{\eni[2]{2},\eni[2]{1},\eni[2]{1}}}
must also belong to the same clone, i.e.\ to \m{\ncent[2]{\set{T}}}.
However, the latter set has been completely described in
Lemma~\ref{lem:Tstar2} above, it contains precisely~65 functions.
Thus, the behaviour of~$g$ on tuples of the form~\m{(x,x,y)} has to
coincide with one of these 65 functions, likewise, the results on tuples
of the form~\m{(x,y,x)} and of the form~\m{(y,x,x)} are determined by
one of these functions, respectively. Moreover, on the three tuples of
the form \m{(x,x,x)}, the three binary operations from
\m{\ncent[2]{\set{T}}} have to prescribe non\dash{}contradictory values.
Therefore, except for the six
tuples that are permutations of\/ \m{(0,1,2)} the values of \m{g} are
determined by one of at most \m{65^3} choices. Altogether no more than
\m{65^3\cdot 3^6 = 200\,201\,625} ternary functions have to be considered.
\par
This can be done by a computer, resulting in a list\footnote{%
This list can be computed using the function
\texttt{findallternaries\_optimised()} from the ancillary file
\texttt{commutationTs.cpp}, and it is given in the file
\texttt{Tcent3\_sorted.txt}.}
of exactly $1\,048\,578$~functions belonging to~\m{\ncent[3]{\set{T}}}.
Again for each of these ternary operations it is readily verified by a
computer that they commute\footnote{%
This verification can be carried out using the function
\texttt{readTcent3("Tcent3\_sorted.txt")} from the ancillary file
\texttt{commutationTs.cpp} and confirms once more the concluding
sentence in the proof of Lemma~\ref{lem:unique-bin-func}.}
with the binary operation~\m{f} given in
Lemma~\ref{lem:unique-bin-func}. Consequently, by a complete case
distinction, we have indeed that
\m{f\in\set{T}^{*(3)*}}.
Together with Lemma~\ref{lem:Tclone2}, this proves
\m{f\in\set{T}^{*(3)*(2)}\setminus\genClone[2]{\set{T}}}
for \m{A=\set{0,1,2}}.
\end{remark}

It is not a suitable strategy to continue indefinitely with
individual verifications that the unique binary operation~$f$ from
Lemma~\ref{lem:unique-bin-func} belongs to more and more accurate upper
approximations \m{\set{T}^{*(\ell)*}},
\m{\ell\rightarrow\infty}, of \m{\bicent{\set{T}}}. Instead we need a
more creative Galois theoretic argument to be sure that
\m{f\in\bicent{\set{T}}}. This confirmation is given in the following
lemma in the form of a primitive positive definition. As it turns out,
the argument used there for \m{k=\abs{A}=3} and the definition of~$f$
from Lemma~\ref{lem:unique-bin-func} can then be generalised to any
\m{k\geq 3}, see Theorem~\ref{thm:pp-defining-separating-function}.
However, we think it is instructive to first show where the idea for
the theorem originates from.

\begin{lemma}\label{lem:pp-formula}
The binary function \m{f \in\set{T}^{*(2)*(2)}} defined
in Lemma~\ref{lem:unique-bin-func} indeed belongs to~\m{\bicent{\set{T}}} for
its graph is definable by a primitive positive formula%
\footnote{The correctness of this formula has been checked with the
Z3-solver~\cite{deMouraBjoernerZ3EfficientSMTsolver,Z3}, see
the script \texttt{checkformulaforbinfunc.z3} available as an ancillary file.}
over \m{A=\set{0,1,2}} involving only the graph of\/~$T$:
\begin{align*}
\bigl\{(&x_2,x_3,x_5)\in A^3 \mathrel{\big\vert} f(x_2,x_3)=x_5\bigr\}\\
&=\lset{(x_2,x_3,x_5)\in A^3}{\exists x_1,x_4\in A\colon
\begin{aligned}[c]
            T(x_1,x_2,x_3,x_4) &= x_5 \land{}\\
            (x_2,x_3,x_2,x_3,x_1,x_2,x_4,x_3) &\in \ker(T) \land{}\\
            (x_3,x_2,x_3,x_2,x_1,x_3,x_4,x_2) &\in \ker(T)
                            \end{aligned}}\\
&=\lset{(x_2,x_3,x_5)\in A^3}{\exists x_1,x_4,u,v\in A\colon
                          \begin{aligned}[c]
                            T(x_1,x_2,x_3,x_4) &= x_5 \land{}\\
                            T(x_2,x_3,x_2,x_3) &= u \land{}\\
                            T(x_1,x_2,x_4,x_3) &= u \land{}\\
                            T(x_3,x_2,x_3,x_2) &= v \land{}\\
                            T(x_1,x_3,x_4,x_2) &= v
                          \end{aligned}}
\end{align*}
\end{lemma}
\begin{proof}
The idea how to construct the graph of~$f$ is by considering the full graph
of~$T$, that is, the relation
\[\lset{(x_1,x_2,x_3,x_4,x_5)\in A^5}{T(x_1,x_2,x_3,x_4) = x_5},\]
and to project it to the second, third and fifth coordinate. This is
motivated by the fact that~$T$ sends only two arguments, $(1,1,2,2)$ and
$(1,2,1,2)$, to one and every other quadruple to zero, and the middle two
components of the two mentioned quadruples coincide with those pairs that
are mapped to one by~$f$. Of course, such a projection will not result in a
function graph, but it almost does. The pairs $(1,2)$ and $(2,1)$ will be
assigned two values each: the value one (as desired for~$f$) and an
erroneous value zero caused by some other quadruples $(x_1,x_2,x_3,x_4)$
with the same middle component $(1,2)$ or $(2,1)$. Hence, the goal is to
remove those quadruples from the relation before projecting. There are 16
disturbing argument tuples in the graph of~$T$ altogether:
\[ \lset{(u,a,b,v)}{\set{a,b}=\set{1,2}, u,v\in A}\setminus\set{(1,1,2,2),(1,2,1,2)}.\]
They need to be removed by imposing additional conditions that have to be
satisfied by the quadruples $(1,2,1,2)$ and $(1,1,2,2)$ since we have to
ensure that these are kept in the relation.

It turns out that this is
possible by imposing just two additional requirements involving the
kernel of~$T$. The kernel is an equivalence relation on quadruples that
we interpret as an octonary relation on~$A$, and it partitions $A^4$ into
two classes: $\set{(1,2,1,2),(1,1,2,2)}$ and the complement~$B$ of this
set in~$A^4$. In particular~$B$ includes all tuples containing a zero or
three ones or three twos or a two in the first position or a one in the
last position.
Using this observation it is easy to verify that the following two sets
jointly (i.e.\ their intersection) exclude all 16 undesired quadruples. So
these two sets represent the restrictions that we are going to
apply to the graph of~$T$:
\begin{align*}
\lset{(x_1,\dots,x_4)\in A^4}{T(x_2,x_3,x_2,x_3) = T(x_1,x_2,x_4,x_3)}
&=A^4\setminus\set{
\begin{array}{@{}*{8}{c@{\,}}c@{}}
0&0&0&1&1&1&2&2&2\\
1&1&1&1&1&2&1&1&1\\
2&2&2&2&2&2&2&2&2\\
0&1&2&0&1&1&0&1&2
\end{array}}\\
\lset{(x_1,\dots,x_4)\in A^4}{T(x_3,x_2,x_3,x_2) = T(x_1,x_3,x_4,x_2)}
&=A^4\setminus\set{
\begin{array}{@{}*{8}{c@{\,}}c@{}}
0&0&0&1&1&1&2&2&2\\
2&2&2&2&2&2&2&2&2\\
1&1&1&1&1&2&1&1&1\\
0&1&2&0&1&1&0&1&2
\end{array}}
\end{align*}
Both sets also exclude the tuple $(1,2,2,1)$, but this is not harmful, as there
are sufficiently many other quadruples left having $(2,2)$ as their middle
component, for example $(0,2,2,0)$.
\end{proof}

As the arity of~$T$ is \m{(k-1)^2} where \m{k=\abs{A}}, it is perhaps
helpful to arrange the arguments of~$T$ in a
\nbdd{((k-1)\times(k-1))}square. Expressing the primitive positive
formula from Lemma~\ref{lem:pp-formula} using such
\nbdd{(2\times2)}squares then yields
\[
\exists\,x_1,x_4\in \set{0,1,2}\colon
T\apply{\begin{smallmatrix}x_1&x_2\\x_3&x_4\end{smallmatrix}} = x_5
\land
T\apply{\begin{smallmatrix}x_2&x_3\\x_2&x_3\end{smallmatrix}} =
T\apply{\begin{smallmatrix}x_1&x_2\\x_4&x_3\end{smallmatrix}}
\land
T\apply{\begin{smallmatrix}x_3&x_2\\x_3&x_2\end{smallmatrix}} =
T\apply{\begin{smallmatrix}x_1&x_3\\x_4&x_2\end{smallmatrix}}.
\]
This kind of interpretation is key for the understanding of the
following main result.

\begin{theorem}\label{thm:pp-defining-separating-function}
Let $A=\set{0,\dotsc,k-1}$ where $k\geq 3$ and put $n=k-1$.
Let the function $f\colon A^n\to A$ be defined by
\[f(\bfa{x}) = \begin{cases}
  1 &\text{if } \bfa{x}\in\set{\asc,\desc},\\
  0 &\text{else},
  \end{cases}\]
where $\asc=(1,\dotsc,n)$ and $\desc=(n,\dotsc,1)$.
The graph of~$f$ can be defined by a primitive positive formula using the
graph of\/~$T$ as follows:
%
\begin{align*}
\bigl\{(&\downdiag,y)\in A^k \mathrel{\big\vert}
f(\downdiag)=y\bigr\}\\
&=\lset{(\downdiag,y)\in A^k}{\apply{\exists
x_{ij}\in A}_{\substack{1\leq i,j\leq n\\i+j\neq k}}\colon
\begin{aligned}
T(\rightarrow_1,\rightarrow_2,\dotsc,\rightarrow_n) &= y\\
T(\downdiag,\downdiag,\dotsc,\downdiag)
&=T(\rightarrow_1,\leftarrow_2,\dotsc,\leftarrow_n)\\
T(\updiag,\updiag,\dotsc,\updiag)
&=T(\downarrow_1,\uparrow_2,\dotsc,\uparrow_n)
\end{aligned}}\\
&=\lset{(\downdiag,y)\in A^k}{
  \apply{\exists x_{ij}\in A}_{\substack{1\leq i,j\leq n\\i+j\neq k}}\,
  \exists u,v\in A\colon
    \begin{aligned}[c]
      T(\rightarrow_1,\rightarrow_2,\dotsc,\rightarrow_n) &= y \land{}\\
      T(\downdiag,\downdiag,\dotsc,\downdiag)             &= u \land{}\\
      T(\rightarrow_1,\leftarrow_2,\dotsc,\leftarrow_n)   &= u \land{}\\
      T(\updiag,\updiag,\dotsc,\updiag)                   &= v \land{}\\
      T(\downarrow_1,\uparrow_2,\dotsc,\uparrow_n)        &= v
    \end{aligned}},
\end{align*}
where the arrows represent the following sequences of variables for
$1\leq i\leq n$:
\begin{align*}
\downdiag &= x_{1,n},x_{2,n-1},\dotsc,x_{n-1,2},x_{n,1}\\
\updiag &= x_{n,1},x_{n-1,2},\dotsc,x_{2,n-1},x_{1,n}\\
\rightarrow_i &= x_{i,1},\dotsc,x_{i,n}\\
\leftarrow_i &= x_{i,n},\dotsc,x_{i,1}\\
\downarrow_i &= x_{1,i},\dotsc,x_{n,i}\\
\uparrow_i &= x_{n,i},\dotsc,x_{1,i}
\end{align*}
\end{theorem}
\begin{proof}
We imagine the $n^2$ variables of~$T$ arranged in a square as follows
\[\Sqre = \begin{matrix}x_{1,1},\dotsc,x_{1,n}\\
\vdots\\
x_{n,1},\dotsc,x_{n,n}\end{matrix},
\]
which we feed row-wise into~$T$, that is, as a notational convention we
identify $\Sqre$ with $\rightarrow_1,\dotsc,\rightarrow_n$ and thus
stipulate
$T(\Sqre)\defeq
T(\rightarrow_1,\dotsc,\rightarrow_n) = T(x_{1,1},\dotsc,x_{n,n})$.
Reversing this line of thought, we can as well start with some
square~$\Sqre$ of
variables, feed its elements into~$f$ in some order (indicated, for
instance, by certain arrows) and then interpret this sequence of
variables as rows of a new square. For example, given~$\Sqre$, the value
$T(\downarrow_1,\dotsc,\downarrow_n)$ is the result of~$T$ applied to a
square whose rows are the columns of~$\Sqre$; so we apply~$T$ to the
transposed~$\Sqre$. Subsequently, we shall often consider sequences as
squares where the rows are connected to the ordering of the given
sequence and the meaning of columns, diagonals etc.\ is tied to this
particular square interpretation.
\par
Two squares play a special role for~$T$, namely those where~$T$
outputs~$1$. First, we have $T(p_1) = 1$ where $p_1$ is given by
${\rightarrow_i} = (i,\dotsc,i)$ for all $1\leq i\leq n$ (that is,
${\downarrow_i}=\asc$ for all $1\leq i\leq n$ and also $\downdiag=\asc$).
Second we have $T(\asc,\dotsc,\asc) = 1$, and we denote the square all
of whose rows $\rightarrow_i$ are $\asc$ by~$p_2$ (this means
${\downarrow_i} = (i,\dotsc,i)$ for all $1\leq i\leq n$ and
$\downdiag=\desc$).
\par
With the square interpretation in mind we form the set
\[\theta = \lset{(\Sqre,y)\in A^{n^2+1}}{
\begin{aligned}
T(\Sqre) &= y\\
T(\downdiag,\downdiag,\dotsc,\downdiag)
&=T(\rightarrow_1,\leftarrow_2,\dotsc,\leftarrow_n)\\
T(\updiag,\updiag,\dotsc,\updiag)
&=T(\downarrow_1,\uparrow_2,\dotsc,\uparrow_n)
\end{aligned}}\]
and then project it to the diagonal $\downdiag$ and the last
coordinate~$y$,
representing the image value of~$T$. To show that this projection
coincides with the graph of~$f$, we shall prove the following
statements:
\begin{enumerate}[(i)]
\item\label{item:p1}
      For every $(\Sqre,y)\in \theta$ where $\downdiag = \asc$, it
      follows $y=1$. This means that $\downdiag=\asc$ implies
      $\Sqre=p_1$.
\item\label{item:p2}
      For every $(\Sqre,y)\in \theta$ where $\downdiag = \desc$, it
      follows $y=1$. This means that $\downdiag=\desc$ implies that
      $\Sqre=p_2$.
\item\label{item:whole-graph-present}
      For every $\bfa{x}\in A^n\setminus\set{\asc,\desc}$ there is
      some~$\Sqre$ such that $(\Sqre,0)\in\theta$ and
      $\downdiag = \bfa{x}$. Moreover, $(p_1,1),(p_2,1)\in\theta$.
\end{enumerate}
Now, if $(\Sqre,y)\in\theta$ then $y\in\im(T) = \set{0,1}$. If $y=1$,
then $\Sqre=p_1$ or $\Sqre=p_2$, whence $\downdiag=\asc$ or
$\downdiag=\desc$ and both $(\asc,1),(\desc,1)\in\graph{f}$. If $y=0$,
then $\Sqre \neq p_1$, so statement~\eqref{item:p1} yields
$\downdiag\neq \asc$; similarly, $\Sqre\neq p_2$ and so
$\downdiag\neq\desc$ by statement~\eqref{item:p2}. Hence in each case we
have $(\downdiag,y)\in\graph{f}$ which shows that the projection
of~$\theta$ is a subset of the graph of~$f$. Conversely,
statement~\eqref{item:whole-graph-present} shows that the full graph
of~$f$ is obtainable as a projection of~$\theta$.
\par
We proceed with the proof of the three statements.
\begin{enumerate}[(i)]
\item If $(\Sqre,y)\in\theta$ and $\downdiag=\asc$, then
      $1= T(\asc,\asc,\dotsc,\asc) =
      T(\rightarrow_1,\leftarrow_2,\dotsc,\leftarrow_n)$.
      This means $(\rightarrow_1,\leftarrow_2,\dotsc,\leftarrow_n)
      \in\set{p_1,p_2}$. Because $\downdiag=\asc$, $x_{1n}=1$ and
      $x_{n1}=n$, so the \nbdd{n}th column of
      $(\rightarrow_1,\leftarrow_2,\dotsc,\leftarrow_n)$ is not constant
      and hence the latter cannot be equal to~$p_2$. Thus it is~$p_1$ and
      therefore also $\Sqre = p_1$.
\item If $(\Sqre,y)\in\theta$ and $\downdiag=\desc$, then reading
      backwards we have $\updiag=\asc$, and therefore
      \m{1= T(\asc,\asc,\dotsc,\asc) =
      T(\downarrow_1,\uparrow_2,\dotsc,\uparrow_n)}, whence
      \m{(\downarrow_1,\uparrow_2,\dotsc,\uparrow_n)\in\set{p_1,p_2}}.
      As $\downdiag=\desc$, we have $x_{1,n}=n$ and $x_{n,1}=1$, so
      the \nbdd{n}th column of
      \m{(\downarrow_1,\uparrow_2,\dotsc,\uparrow_n)} is not constant
      (recall that, by our convention, these tuples are fed as rows
      into~$T$). This means that
      \m{(\downarrow_1,\uparrow_2,\dotsc,\uparrow_n)} must have constant
      rows (be equal to~$p_1$), so \m{{\downarrow_1} = (1,\dotsc,1)}, and
      \m{{\downarrow_i} = (i,\dotsc,i)} for $2\leq i\leq n$. This means
      that~$\Sqre$ has constant columns with values $1,\dotsc,n$, which
      means that $\Sqre=p_2$.
\item First we check that $(p_1,1)\in\theta$. Clearly, $T(p_1)=1$. For
      $p_1$ we have $\downdiag=\asc$ and $\updiag=\desc$, so
      $T(\asc,\dotsc,\asc) = 1 = T(p_1) =
      T(\rightarrow_1,\leftarrow_2,\dotsc,\leftarrow_n)$ holds as $p_1$
      has constant rows, and
      $T(\desc,\dotsc,\desc) = 0=
      T(\asc,\desc,\dotsc,\desc)=
      T(\downarrow_1,\uparrow_2,\dotsc,\uparrow_n)$ is true, as
      well.
      \par
      Next we verify that $(p_2,1)\in\theta$. Again, $T(p_2)=1$.
      This time we have $\downdiag=\desc$ and $\updiag=\asc$, so
      \m{T(\desc,\dotsc,\desc)=0 = T(\asc,\desc,\dotsc,\desc)
      =T(\rightarrow_1,\leftarrow_2,\dotsc,\leftarrow_3)}. Furthermore,
      \m{T(\asc,\dotsc,\asc)=1 =T(p_1)
      =T(\downarrow_1,\uparrow_2,\dotsc,\uparrow_n)} because the columns
      of~$p_2$ have constant values $1,\dotsc,n$.
      \par
      Finally, consider some $\bfa{x}\in A^n\setminus\set{\asc,\desc}$
      and $\Sqre$ with $\downdiag=\bfa{x}$ and having zeros everywhere
      else.
      All rows of
      $(\downdiag,\dotsc,\downdiag)$ and of $(\updiag,\dotsc,\updiag)$
      are identical, so none of these two squares is~$p_1$. If one of
      these were~$p_2$, then $\bfa{x}=\downdiag=\asc$ or $\updiag=\asc$,
      which would mean $\bfa{x}=\downdiag=\desc$. Both options are
      excluded by the choice of~$\bfa{x}$. Since neither of these two
      squares is $p_1$ or~$p_2$, we have
      $T(\downdiag,\dotsc,\downdiag)=0=T(\updiag,\dotsc,\updiag)$.
      As $\Sqre$ has zeros outside the \nbdd{\downdiag}diagonal, it
      follows that also
      $(\rightarrow_1,\leftarrow_2,\dotsc,\leftarrow_n)$ and
      $(\downarrow_1,\uparrow_2,\dotsc,\uparrow_n)$ have zeros somewhere
      and are hence mapped to zero by~$T$. Thus~$\Sqre$ satisfies the
      two conditions regarding the kernel of~$T$. As~$\Sqre$ contains
      zeros, we also have $T(\Sqre) = 0 = y$, concluding the argument.
      \qedhere
\end{enumerate}
\end{proof}

As a corollary we obtain that the example algebras $\algwops{A}{T}$
constructed by Snow in~\cite{SnowGeneratingPrimitivePositiveClones} do
not generate centraliser clones as term operations and are thus no
counterexample to the Burris-Willard conjecture or to
Dani\v{l}\v{c}enko's results.

\begin{corollary}\label{cor:f-separating-clones}
For every carrier~$A$ of cardinality $k\geq 3$ the \nbdd{(k-1)}ary
function~$f$ defined in
Theorem~\ref{thm:pp-defining-separating-function} satisfies
$f\in \bicent{\set{T}}\setminus\genClone{\set{T}}$.
\end{corollary}
\begin{proof}
By Theorem~\ref{thm:pp-defining-separating-function} we have
$f\in\bicent{\set{T}}$; since~$f$ is \nbdd{(k-1)}ary, it cannot belong to the
clone generated by~$T$ as it is maps two distinct tuples to one and is
not a projection (cf.\ Lemma~\ref{lem:Tclone2-general}).
\end{proof}

\section{Some computational remarks}\label{sect:computations}
We conclude with a few comments on computational aspects related to
verifying that for \m{A=\set{0,1,2}}, the simplest case in question,
the binary function
\m{f\in\set{T}^{*(2)*(2)}\setminus\genClone[2]{\set{T}}}
found in Lemma~\ref{lem:unique-bin-func} actually belongs to
\m{\bicent{\set{T}}}.
\par

The first possibility is based on trusting the classification results
shown by Da\-ni\v{l}\-\v{c}en\-ko in~\cite[Theorems~4, 5,
pp.~103, 105]{Danilcenko1979-thesis}. Using the equivalence of
statements~\eqref{item:cent-leq-n} and~\eqref{item:cent-n-cent-bicent}
in Proposition~\ref{prop:char-cdeg}, these theorems imply that
\m{\bicent{\set{T}}=\set{T}^{*(3)*}}, which
contains~$f$ by the calculations described in
Remark~\ref{rem:f-in-cent-3-cent-T}. Believing in Dani\v{l}\v{c}enko's
thesis obviously does not render
Theorem~\ref{thm:pp-defining-separating-function} obsolete, as the
latter also covers the cases where \m{\abs{A}>3}.
\par

The second option we would like to discuss is whether it is feasible to
compute a primitive positive formula over~\m{\graph{T}} that allows to
define~\m{\graph{f}}. The formula shown in
Lemma~\ref{lem:pp-formula} uses five \nbdd{\graph{T}}atoms and
four existentially quantified variables. Of course, these bounds are not
known beforehand, and even if they were, simply trying to produce all
formul\ae{} with $\ell=1,2,3,\dots$ \nbdd{\graph{T}}atoms and trying to
find a \nbdd{3}variable projection that gives~$\graph{f}$ becomes
unwieldy very quickly. Indeed, before even dealing with projections,
there are \m{\kappa^\kappa} possible variable substitutions, where
\m{\kappa\defeq\ell\cdot\arity(\graph{T})} for a primitive positive
formula with $\ell$~atoms of type~$\graph{T}$ and at most~$\kappa$
variables. More concretely, to get the formula from
Lemma~\ref{lem:pp-formula}, we would have \m{\ell=5} and
\m{\arity(\graph{T})=5}, so \m{\kappa=25}, and
\m{25^{25}\approx 10^{35}} substitutions are currently too many to
check in a reasonable amount of time.
\par

However, if \m{f\in\bicent{\set{T}}}, then
\m{\graph{f}\in\Inv{\Pol{\set{\graph{T}}}}}, and there is a more
systematic method to compute a primitive positive formula for a
relation \m{\rho_0\in\Inv{F}} on a finite set~$A$ where
\m{F=\Pol{\set{\rho_1,\dotsc,\rho_t}}}, \m{t\in\N}.
It comes from an algorithm to compute~$\Fn{F}$ interpreted as
a relation~\m{\Gamma_F(\chi_n)} of arity~$\abs{A}^n$, which is given
in~\cite[4.2.5., p.~100 et seq.]{PoeKal}, combined with the proof of
the second part of the main theorem on the $\PolOp\text{-}\InvOp$
Galois connection, showing that any
\m{\rho_0\in\Inv{\Pol{Q}}} belongs to the relational clone generated
by~$Q$, as it is primitive positively definable from
\m{\Gamma_{F}(\chi_n)} where \m{F=\Pol{Q}} and \m{n=\abs{\rho_0}}
(cf.~\cite[1.2.2.~Lemma, p.~53 et seq.]{PoeKal}).
\par

The following is slightly more general than what is described
in~\cite[4.2.5]{PoeKal} for we can deal with finitely many describing
relations \m{\rho_1,\dotsc,\rho_t} for the polymorphism clone~$F$,
while only one is used in~\cite{PoeKal}. Taking
\m{Q=\set{\rho_1\times\dotsm\times\rho_t}} as a singleton
in~\cite{PoeKal} is inefficient from a computational point of view, so
we give a proof of this not very original modification. Additionally,
we allow for a generating system~$\gamma_0$ of the relation~$\rho_0$
for which a formula is sought (although this is somehow implicit
in~\cite[4.2.5]{PoeKal} as \m{\Gamma_{F}(\chi_n)} is generated by the
\nbdd{n}element subrelation~$\chi_n$).

\begin{proposition}\label{prop:gen-rel-clone}
Assume \m{Q\defeq\lset{\rho_\ell}{1\leq \ell\leq t}}, \m{F\defeq \Pol{Q}},
\m{\rho_0\in\Inv{F}} where \m{\rho_\ell\subs A^{m_\ell}} for
\m{0\leq \ell\leq t}. Let \m{\gamma_0\subs\rho_0} with \m{n\defeq
\abs{\gamma_0}} be a generating system of~$\rho_0$, that is,
\m{\rho_0=\gapply{\gamma_0}_{\algwops{A}{F}^{m_0}}}.
There is \m{m'\leq m_0} and \m{\gamma\subs\rho\in\Inv[m']{F}} and
\m{\alpha\colon m_0\to m'} such that
\m{\gamma_0=\lset{x\circ\alpha}{x\in\gamma}} where~$\gamma$ does not have
any duplicate coordinates. If we imagine the tuples in~\m{\gamma_0}
written as columns of an \nbdd{(m_0\times n)}matrix, then the distinct
rows of this matrix are precisely the rows of the matrix whose columns
form the tuples of~$\gamma$. Some of these $m'$~rows will be found as rows of a
relation \m{\mu\subs A^L} with \m{\abs{\mu}=n} defined below. For
notational simplicity we choose~$\alpha$ such that the rows with
indices \m{1,\dotsc,m} have this property and put \m{p\defeq
m'-m\geq0}.
\par
The matrix representation of the relation~$\mu$ has~$n$ columns
(tuples) and~$L$ rows \m{\apply{\bfa{z}_i}_{0\leq i<L}} where \m{L=\sum_{\ell=1}^t s_\ell^n\cdot m_\ell}
with \m{s_\ell=\abs{\rho_\ell}}. Let the columns of~$\mu$ arise by
stacking on top of each other all possible submatrices of~$\rho_1$
with~$n$ columns, followed by all possible submatrices of~$\rho_2$,
and so forth, finishing with all submatrices obtained by choosing~$n$
of the~$s_t$ columns of~$\rho_t$. Thus
\m{\mu\subs\pi\defeq\rho_1^{s_1^n}\times\dotsm\times\rho_t^{s_t^n}}.
Define the kernel relation \m{\epsilon\defeq \lset{(i,j)\in
L^2}{\bfa{z}_i = \bfa{z}_j}} and identify variables in~\m{\pi}
accordingly with \m{\delta_\epsilon = \lset{x\in A^L}{\forall\,
(i,j)\in\epsilon\colon x_i=x_j}}. This gives
\m{\sigma\defeq\pi\cap\delta_{\epsilon}} having the same row kernel as~$\mu$.
By finding the first~$m$ rows of~$\gamma$ among the rows of~$\mu$, we
find a projection~$\pr$ to an \nbdd{m}element set of indices such that
\m{\gamma\subs\pr(\mu)\times A^p}. It follows that
\m{\rho=\pr(\sigma)\times A^p
       =\pr\apply{\apply{\rho_1^{s_1^n}\times
       \dotsm\times \rho_t^{s_t^n}}\cap\delta_\epsilon}\times A^p}.
\end{proposition}
\begin{proof}
To show that \m{\rho\subs\pr(\sigma)\times A^p} we note that
\m{\pr(\sigma)\times A^p\in \Inv{F}} since for every \m{1\leq\ell\leq t}
we have \m{\rho_\ell\in Q\subs\Inv{F}}. Moreover, as~$\rho$ is a
projection of~$\rho_0$ in the same way as~$\gamma$ is a projection
of~$\gamma_0$, and since
\m{\rho_0=\gapply{\gamma_0}_{\algwops{A}{F}^{m_0}}}, we have
\m{\rho=\gapply{\gamma}_{\algwops{A}{F}^{m'}}}.
Due to \m{\gamma\subs\pr(\mu)\times A^p \subs \pr(\sigma)\times A^p},
the generating set~$\gamma$ is a subset of the invariant
\m{\pr(\sigma)\times A^p}, and so is the generated invariant~$\rho$.
\par
For the converse inclusion we take a tuple $x\in\sigma$ and some
\m{a\in A^p} and denote by~$y$ the tuple obtained from~$x$ by
projection to the~$m$ indices that relate the first~$m$ rows
\m{\bfa{v}_1,\dotsc,\bfa{v}_m} of~$\gamma$ to a certain section of~$x$.
Since \m{x\in\sigma\subs\delta_\epsilon}, this tuple defines a function
\m{f_x\colon B\to A} where
\m{B\defeq\lset{\bfa{z}_i}{0\leq i<L}\sups
   C\defeq\lset{\bfa{v}_i}{1\leq i\leq m}} by finding for any
\m{\bfa{z}\in B} some \m{0\leq i<L} such that
\m{\bfa{z}=\bfa{z}_i} and letting \m{f_x(\bfa{z}) \defeq x_i} (the
choice of \m{i<L} is inconsequential as \m{x\in\delta_\epsilon}).
The remaining rows \m{\bfa{v}_i} with \m{m<i\leq m+p} do not belong
to~$B$ by the choice of~$m$ and~$p$. This means, $(y,a)$ defines a
function \m{f_{y,a}\colon\lset{\bfa{v}_i}{1\leq i\leq m'}\to A}, where
\m{f_x\restriction_C=f_{y,a}\restriction_C}. Moreover, it is possible
to extend~$f_{y,a}$ to a globally defined function \m{f\colon A^n\to A}
such that \m{f\restriction_{\set{\bfa{v}_i\mid 1\leq i\leq m'}}=f_{y,a}}
and \m{f\restriction_B=f_x} without contradictory value assignments.
We pick one particular such~$f$, no matter which one, and we show below
that~$f\in F=\Pol{Q}$. By the hypothesis of the proposition, $\rho$
belongs to~$\Inv{F}$, so~$f$ preserves~$\rho$. Thus, applying~$f$ to
the tuples in~$\gamma$, gives
\m{(y,a)=\apply{f_{y,a}(\bfa{v}_i)}_{1\leq i\leq m'}
        =\apply{f(\bfa{v}_i)}_{1\leq i\leq m'}\in\rho} as needed.
\par
It remains to argue that~$f\in\Pol{Q}$. Hence, take any
\m{1\leq \ell\leq t} and any matrix of $n$~columns taken
from~$\rho_\ell$.  By the construction of~$\mu$ there are
$m_\ell$~consecutive indices \m{0\leq i,i+1,\dotsc,i+m_\ell-1<L} such
that the rows of this matrix are
\m{\bfa{z}_i,\dotsc,\bfa{z}_{i+m_\ell-1}}.
Now \m{\apply{f(\bfa{z}_{i+\nu})}_{0\leq \nu<m_\ell}=
\apply{f_x(\bfa{z}_{i+\nu})}_{0\leq \nu<m_\ell}}, and this tuple is
in~$\rho_\ell$ because~$f_x$ is defined via
\m{x\in\sigma=\pi\cap\delta_{\epsilon}}.
\end{proof}

The expression \m{\rho=\pr\apply{\apply{\rho_1^{s_1^n}\times
       \dotsm\times \rho_t^{s_t^n}}\cap\delta_\epsilon}\times A^p}
in Proposition~\ref{prop:gen-rel-clone} gives a primitive positive
definition of~$\rho$ in terms of~$\rho_1,\dotsc,\rho_t$. Duplicating
variables as indicated by~$\alpha$, one can then give a primitive
positive formula for the original relation~$\rho_0$.
The inclusion \m{\rho\subs\pr\apply{\apply{\rho_1^{s_1^n}\times
       \dotsm\times \rho_t^{s_t^n}}\cap\delta_\epsilon}\times A^p}
holds in any case, regardless of the assumption that
\m{\rho\in\Inv{F}}. It can be seen from the proof of
Proposition~\ref{prop:gen-rel-clone} that the latter condition is only
needed for the opposite inclusion.
\par
That is, the formula computed by the following algorithm will always be
satisfied by all tuples from~$\rho$ (or~$\rho_0$), but if the
containment
\m{\rho\in\Inv{F}} is only suspected but not known in advance, then one
needs to check afterwards that the tuples satisfying the
generated primitive positive formula really belong to~$\rho$ (or
to~$\rho_0$, respectively).

\begin{algo}\label{alg:gen-rel-clone}
Compute a primitive positive definition\footnote{%
An implementation is available in the file \texttt{ppdefinitions.cpp},
which can be compiled using \texttt{compile.sh}, resulting in an
executable \texttt{getppformula}. This executable expects a file
\texttt{input.txt}, the formatting of which is explained in
\texttt{input\_template.txt}, which can also be used as
\texttt{input.txt}. After a successful run the programme will
produce files \texttt{ppoutput.out}, an ascii text file containing the
computed primitive positive formula, and \texttt{checkppoutput.z3}, a
script to verify the correctness of the formula using the Z3 theorem
prover~\cite{deMouraBjoernerZ3EfficientSMTsolver,Z3}.
\par
There are two caveats with the implementation added as ancillary file
to this submission: first, the initial preprocessing step
turning~$\gamma_0$ into~$\gamma$ has not been implemented. Hence,
\texttt{ppdefinitions.cpp} expects a goal relation~$\gamma$
(relation~\texttt{S} in \texttt{input.txt}) without duplicate
coordinates. If $\gamma_0$ has duplicate rows, the initial massaging
and the final adjustment of the formula by duplicating the respective
variables has to be done by hand. Second, it is possible to use a
proper generating set $\gamma_0\subs\rho_0$ in the input (provided it
does not contain duplicate rows), but then in the output file
\texttt{checkppoutput.z3} the goal relation~\texttt{S} has to be
completed manually with all tuples from~$\rho_0$, since the closure
\m{\gapply{\gamma_0}_{\algwops{A}{F}^{m_0}}} is not computed.}
\newline
(Pseudocode is given on page~\pageref{code:compute-ppdefinition},
line numbers in the description refer to this code.)
\begin{description}
\item[Input]
  finitary relations \m{\rho_1\subs A^{m_1},\dotsc,\rho_t\subs A^{m_t}}
  defining \m{F\defeq \Pol{Q}} where
  \m{Q=\set{\rho_\ell\mid 1\leq \ell\leq t}}
  \par
  a generating system \m{\gamma_0\subs A^{m_0}} for a relation
  \m{\rho_0=\gapply{\gamma_0}_{\algwops{A}{F}^{m_0}}\in \Inv{F}}
\item[Output]
  a primitive positive formula describing \m{\rho_0} in terms of
  \m{\rho_1,\dotsc,\rho_t}
\item[Description]
  We assume that \m{\gamma_0=\set{r_1,\dotsc,r_n}}, the tuples of which
  we represent as a matrix with columns \m{r_1,\dotsc,r_n} and rows
  \m{v_1,\dotsc,v_{m_0}}.
  We first define a map
  \m{\alpha\colon\set{1,\dotsc,m_0}\to\set{1,\dotsc,m'}} to a
  transversal of the equivalence relation
  \m{\lset{(i,j)\in\set{1,\dotsc,m_0}^2}{v_i=v_j}} (lines~1--9).
  For this we iterate over all rows, and, if~$v_j$ has been seen
  previously among \m{v_1,\dotsc,v_{j-1}}, we assign to \m{\alpha(j)}
  the same index \m{\iota(v_j)} as previously, and if~$v_j$ is a fresh
  row, we assign to~\m{\alpha(j)} the least index \m{i\eqdef\iota(v_j)}
  not used before (lines~4--9). When this is finished,
  \m{\gamma_0=\lset{(x_{\alpha(1)},\dotsc,x_{\alpha(m_0)})}{(x_1,\dotsc,x_{m'})\in\gamma}}
  where \m{\gamma\subs A^m} is a projection of~$\gamma_0$ to its
  distinct rows, and~$m'$ is the last used value of~$i$.
  \par
  Next we iterate over all \m{1\leq \ell\leq t} and for each
  relation~$\rho_\ell$ we iteratively extend the set~$\mathcal{L}_\ell$
  of \nbdd{\rho_\ell}atoms for the final formula, starting from
  \m{\mathcal{L}_\ell=\emptyset} (lines~10--13).
  We iterate over the rows \m{z_1,\dotsc,z_{m_\ell}} of all possible
  matrices with~$n$ columns chosen from~$\rho_\ell$ (lines~14--16). For any of these
  matrices we construct an \nbdd{m_\ell}tuple~$a$ of variable symbols (lines~17--24),
  which will represent a \nbdd{\rho_\ell}atom and will be added
  to~$\mathcal{L}_\ell$ if it is not already present in the list of
  atoms (lines~25--26). The atoms have to be constructed in such a way that any two
  identical rows occurring within all possible matrices get the same
  variable symbol. This ensures that the variable identification
  represented in Proposition~\ref{prop:gen-rel-clone} by intersection
  with~$\delta_\epsilon$ takes place. Moreover, if a row in the
  matrices occurs as a row of~$\gamma$ (or equivalently of~$\gamma_0$),
  then the corresponding variable is not going to be existentially
  quantified, while all others are. This takes care of the projection
  in the formula for~$\rho$ from Proposition~\ref{prop:gen-rel-clone}.
  \par
  In more detail, if a row~$z_j$ with \m{1\leq j\leq m_\ell} has not occurred
  previously (line~17), we have to define its variable
  symbol~$u(z_j)$. If \m{z_j\in\set{v_1,\dotsc,v_{m_0}}}, that is,
  $z_j$ is among the rows of~$\gamma_0$, we use the
  variable \m{u(z_j)\defeq x_{\iota(z_j)}} (lines~19--20).
  Otherwise, the fresh row~$z_j$ needs to be projected away by
  existential quantification, and we use a different symbol
  \m{u(z_j)\defeq y_k} where~$k>0$ is the least previously unused index
  for existentially quantified variables (lines~21--23). Regardless of whether~$z_j$
  is fresh or not, we define the \nbdd{j}th entry of the current
  atom~$a$ as \m{a(j)\defeq u(z_j)} (line~24). Only if the resulting string
  $a=(a(1),\dotsc,a(m_\ell))\notin\mathcal{L}_\ell$, that is, if~$a$ is
  a new atom, it will be added to~$\mathcal{L}_\ell$ (lines~25--26).
  \par
  After all iterations, we state that all variables \m{x_1,\dotsc,x_i}
  occurring in \m{\set{x_{\alpha(1)},\dotsc,x_{\alpha(m_0)}}} come from
  the base set~$A$, we existentially quantify all variables
  \m{y_1,\dotsc,y_k} and write out (line~27) a long conjunction over all
  relations \m{\rho_1,\dotsc,\rho_t} and over all \nbdd{\rho_\ell}atoms
  \m{a\in \mathcal{L}_\ell}
  (cf.\ the direct product in the formula for~$\rho$ in
  Proposition~\ref{prop:gen-rel-clone}).
\end{description}
\begin{algorithm}[htp]
 \SetAlgoVlined
 \SetKwInOut{Input}{Input}
 \SetKwInOut{Output}{Output}
\Input{%
  finitary relations \m{\rho_1\subs A^{m_1},\dotsc,\rho_t\subs A^{m_t}}\\
  \texttt{// defining \m{F\defeq \Pol{Q}} where
  \m{Q=\set{\rho_\ell\mid 1\leq \ell\leq t}}}\\
  generating system $\gamma_0\subs A^{m_0}$ for a relation
  \rlap{$\rho_0=\gapply{\gamma_0}_{\algwops{A}{F}^{m_0}}\in\Inv{F}$}\\
  \tcp{where \m{\gamma_0=\set{r_1,\dotsc,r_n}}, i.e.,
       \m{\abs{\gamma_0}\leq n}
       \newline
       written as a matrix
         $\apply{r_1,\dotsc,r_n} = \apply{\begin{smallmatrix}
         v_1\\\vdots\\v_{m_0}
         \end{smallmatrix}}$
         with rows $v_j\in A^n$}
  }
\Output{a primitive positive presentation of~$\rho_0$ in terms of
  \m{\rho_1,\dotsc,\rho_t}}
\Begin{%
$i\gets 0$
\tcp*[r]{initialise index for distinct rows of~$\gamma_0$}
$D_0\gets\emptyset$
\tcp*[r]{initialise domain of distinct rows of~$\gamma_0$}
\tcp{Define $\iota\colon D_0\to \set{1,\dotsc,\abs{D_0}}$,
            $\alpha\colon \set{1,\dotsc,m_0}\to\set{1,\dotsc,\abs{D_0}}$}
\ForAll{$1\leq j\leq m_0$}{%
  \If{$v_j\notin D_0$}{
    $D_0\gets D_0\cup\set{v_j}$\;
    $i\gets i+1$\;
    $\iota(v_j)\gets i$\;
  }
  $\alpha(j)\gets\iota(v_j)$
}
\tcp{Now \m{D_0=\set{v_1,\dotsc,v_{m_0}}}}
$k\gets 0$
\tcp*[r]{initialise index for $\exists$-quantified variables}
$D\gets \emptyset$
\tcp*[r]{\mbox{initialise domain of distinct rows from submatrices}
         of~$\rho_1,\dotsc,\rho_t$ to define
         $u\colon D\to \set{x_{1},\dotsc,
         x_{\abs{D_0}}}\cup\set{y_1,\dotsc,y_{k}}$}
\ForAll{$1\leq \ell\leq t$}{
\m{\mathcal{L}_\ell\gets \emptyset}
\tcp*[r]{initialise list of atoms pertaining to~$\rho_\ell$}
\ForAll{$c\colon n\to\rho_\ell$}{%
Form a matrix
\m{\apply{c_0,\dotsc,c_{n-1}}=\apply{\begin{smallmatrix}z_1\\\vdots\\z_{m_\ell}\end{smallmatrix}}}
with rows \m{z_j\in A^n}\;%
\tcp{Iterate over its rows and form a possibly new atom~$a$}
\ForAll{$1\leq j\leq m_\ell$}{%
\If(\tcp*[f]{A previously unseen row $z_j$ appears.}){$z_j\notin D$}{
  $D\gets D\cup\set{z_j}$\;
  \eIf(\tcp*[f]{It is a row of~$\gamma_0$.}){$z_j\in D_0$}{%
    $u(z_j)\gets x_{\iota(z_j)}$}
  {
    $k\gets k+1$\;
    $u(z_j)\gets y_{k}$}
}
$a(j)\gets u(z_j)$
\tcp*[r]{extend current atom with the appropriate variable symbol}
}
\If(\tcp*[f]{If it is really new\dots}){$a=(a(1),\dotsc,a(m_\ell))\notin \mathcal{L}_\ell$}{
  $\mathcal{L}_\ell \gets \mathcal{L}_\ell\cup\set{a}$
  \tcp*[r]{\dots add current atom~$a$ to the list.}
}
}}
\Return{String
\m{\rho_0=\Bigl\{\apply{x_{\alpha(1)},\dotsc,x_{\alpha(m_0)}}
\ \Big\vert\ x_1,\dotsc,x_i\in A\land
\exists y_1\dotsm \exists y_k\colon
\bigwedge\limits_{1\leq \ell\leq t}\bigwedge\limits_{a\in\mathcal{L}_\ell}\rho_\ell(a)\Bigr\}}}
}
\NoCaptionOfAlgo
\caption{Compute a primitive positive definition}
\label{code:compute-ppdefinition}
\end{algorithm}
\end{algo}

\begin{example}\label{ex:computing-ternary-f.-from-T.}
In the case discussed in this section, we have \m{A=\set{0,1,2}},
\m{t=1}, \m{Q=\set{\graph{T}}}, \m{\rho_0=\graph{f}}, \m{m_0=3} and
\m{m_1=5}. Moreover, \m{F=\Pol{Q}=\cent{\set{T}}}.
As the size~$s_1$ of~$\graph{T}$ is \m{\abs{A}^4=81}, it is
crucial for the applicability of Algorithm~\ref{alg:gen-rel-clone} to
find a small generating system~$\gamma_0$ of~$\graph{f}$ with respect to
\m{\alg{A}^3} where \m{\alg[\cent{\set{T}}]{A}}. Given
\m{\abs{\gamma_0}=n}, the algorithm has to iterate over \m{s_1^n=81^n}
matrices and thus over \m{m_1\cdot s_1^n = 5\cdot 81^n} rows.
Experiments show that if we blindly took \m{\gamma_0=\graph{f}}, i.e.,
\m{n=\abs{A}^2=9}, the algorithm would need more than eighteen thousand
years to finish, perhaps less by a factor of ten if run on a computer
much faster than the author's.
Fortunately, the number~$n$ can be reduced significantly to a value far
below~$9$.
\par
Indeed, listing the tuples of~$\graph{f}$ as columns, we have
\begin{align*}
\graph{f}&=\set{
\apply{\begin{smallmatrix}0\\0\\0\end{smallmatrix}},
\apply{\begin{smallmatrix}0\\1\\0\end{smallmatrix}},
\apply{\begin{smallmatrix}0\\2\\0\end{smallmatrix}},
\apply{\begin{smallmatrix}1\\0\\0\end{smallmatrix}},
\apply{\begin{smallmatrix}1\\1\\0\end{smallmatrix}},
\apply{\begin{smallmatrix}1\\2\\1\end{smallmatrix}},
\apply{\begin{smallmatrix}2\\0\\0\end{smallmatrix}},
\apply{\begin{smallmatrix}2\\1\\1\end{smallmatrix}},
\apply{\begin{smallmatrix}2\\2\\0\end{smallmatrix}}}\\
&=\gapply{\set{
\apply{\begin{smallmatrix}1\\2\\1\end{smallmatrix}},
\apply{\begin{smallmatrix}2\\1\\1\end{smallmatrix}}
}}_{\alg{A}^3}.
\end{align*}
To see this, we can take advantage of the unary operations
\m{u_{2,a}\in \ncent[1]{\set{T}}} with \m{a\in A}, described in
Corollary~\ref{cor:Tstar1}, and
\m{f_{0,(2,2,2,2)}\in\ncent[2]{\set{T}}} from
Lemma~\ref{lem:Tstar2}. Namely, for \m{a\in\set{0,1,2}}, we have
\begin{align*}
u_{2,a}(1)&=0& u_{2,a}(2)&=a& f_{0,(2,2,2,2)}(1,2)&=2& u_{2,1}(2)&=1\\
u_{2,a}(2)&=a& u_{2,a}(1)&=0& f_{0,(2,2,2,2)}(2,1)&=2& u_{2,1}(2)&=1\\
u_{2,a}(1)&=0,& u_{2,a}(1)&=0,& f_{0,(2,2,2,2)}(1,1)&=0,& u_{2,1}(0)&=0.
\end{align*}
\par
Hence, we can use the \nbdd{2}element generating set
\m{\gamma_0=\set{(1,2,1),(2,1,1)}\subs\graph{f}}
and thus we only have to enumerate \m{5\cdot 81^2 = 32\,805} rows. This
can be done in a fraction of a second\footnote{%
After compilation the programme \texttt{ppdefinitions.cpp} may be run
on \texttt{input\_2generated.txt} copied to \texttt{input.txt}. The
mentioned files can be found in the ancillary directory of this submission.}
and results in a primitive positive formula\footnote{%
Running \texttt{ppdefinitions.cpp} on \texttt{input\_2generated.txt}
(see the ancillary directory) produces the content of \texttt{ppoutput\_2generated.out}
and \texttt{checkppoutput\_2generated.z3}, which both contain the resulting
primitive positive formula (as plain text and in
\texttt{SMT-LIB2.0}-syntax).} with $6$~existentially quantified variables and
$6\,561$~\nbdd{\graph{T}}atoms, the correctness of which can be
verified by a sat\dash{}solver in a few minutes\footnote{%
This can, for example, be done with the Z3 theorem
prover~\cite{deMouraBjoernerZ3EfficientSMTsolver,Z3} using the
ancillary file \texttt{checkppoutput\_2generated.z3}. This file also
contains the computed primitive positive formula for~$\graph{f}$
expressed in the \texttt{SMT-LIB2.0}-format.}.
\end{example}

We conclude that it is possible to computationally find a proof that
the graph of~$f$ is primitive positively definable from~$\graph{T}$ for
\m{A=\set{0,1,2}}. However, the resulting formula is not suitable for a
generalisation to larger carrier sets as the one from
Lemma~\ref{lem:pp-formula} was.

\section*{Acknowledgements}\label{sect:acknowledgements}
The author would like to thank Zarathustra Brady for telling him about
the possibility to include ancillary files with an arXiv preprint.
Moreover, he is grateful to Dmitriy Zhuk for mentioning the usefulness
of the Z3 theorem prover in connection with clones.

\input{referencesBW.tex}%
%
\end{document}

%% file: referencesBW.tex
\def\rus#1{\foreignlanguage{russian}{\bgroup\fontfamily{cmr}\fontencoding{T2A}\selectfont#1\egroup}}
\providecommand{\href}[2]{\texttt{#1}}
\providecommand*{\doi}[1]{\href{http://dx.doi.org/\detokenize{#1}}{\detokenize{#1}}}

%% file: noteBWconjecture.bbl
\begin{thebibliography}{10}

\bibitem{BodnarcukKaluzninKotovRomovGaloisTheoryForPostAlgebras}
V.~G. Bodnarčuk, Lev~Arkaďevič Kalužnin,
  Victor~N. Kotov, and Boris~A. Romov.
\newblock \rus{Теория Галуа для алгебр Поста. I, II}
  [{G}alois theory for {P}ost algebras. {I}, {II}].
\newblock \emph{Kibernetika (Kiev)}, 5(3):1--10; ibid. 1969, no. 5, 1--9, 1969.
\newblock Translated as \doi{10.1007/BF01070906} and \doi{10.1007/BF01267873}.

\bibitem{BurrisWillardFinitelyManyPPClones}
Stanley Burris and Ross Willard.
\newblock Finitely many primitive positive clones.
\newblock \emph{Proc. Amer. Math. Soc.}, 101(3):427--430, November 1987.
\newblock \doi{10.2307/2046382}.

\bibitem{Danilcenko1974ParametricallyClosedClasses3ValuedLogic}
Anna~Fedorovna Daniľčenko.
\newblock \rus{О параметрически замкнутых классах функций трехзначной логики}
  [{O}n parametrically closed classes of functions of three-valued logic].
\newblock In \emph{Proceedings of the 3rd {A}ll-{U}nion {C}onference on
  {M}athematical {L}ogic (Abstracts)}, page~53, Novo\-si\-birsk, 1974.
\newblock In Russian.

\bibitem{Danilcenko1976ParametricallyIndecomposables}
Anna~Fedorovna Daniľčenko.
\newblock \rus{Параметрически неразложимые функции трехзначной логики}
  [{P}arametrically indecomposable functions of three-valued logic].
\newblock In \emph{Proceedings of the 4th {A}ll-{U}nion {C}onference on
  {M}athematical {L}ogic (Abstracts)}, page~38, Kišinev, 1976.
\newblock In Russian.

\bibitem{Danilcenko1977ParametricExpressibility3ValuedLogic}
Anna~Fedorovna Daniľčenko.
\newblock \rus{О параметрической выразимости функций трехзначной логики}
  [{P}arametric expressibility of functions of three-valued logic].
\newblock \emph{Algebra i Logika}, 16(4):397--416, 493, 1977.
\newblock In Russian.

\bibitem{Danilcenko1977ParametricExpressibility3ValuedLogicTranslated}
Anna~Fedorovna Daniľčenko.
\newblock Parametric expressibility of functions of three-valued logic.
\newblock \emph{Algebra and Logic}, 16(4):266--280, July 1977.
\newblock \doi{10.1007/BF01669278}.

\bibitem{Danilcenko1978ParametricallyClosedClasses3ValuedLogic}
Anna~Fedorovna Daniľčenko.
\newblock \rus{Параметрически замкнутые классы функций трехзначной логики}
 [{P}arametrically closed classes of functions of three-valued logic].
\newblock \emph{Bul. Akad. Štiince RSS Moldoven.}, (2):13--20, 93, 1978.
\newblock In Russian.

\bibitem{Danilcenko1979-thesis}
Anna~Fedorovna Daniľčenko.
\newblock \emph{\rus{Вопросы параметрической выразимости функций трехзначной логики}
  [{P}roblems of parametric expressibility of functions in three-valued logic]}.
\newblock PhD thesis, {I}nstitut {M}atematiki s
  {V}yčisliteľnym {C}entrom, {A}kademija {N}auk
  {M}oldavskoj {S}{S}{R}. [{I}nstitute of {M}athematics and {C}omputing Centre
  of the {A}cademy of {S}ciences of the {M}oldovian {S}{S}{R}],
  Kišinev-277028, ul. Akademičeskaja 5 [Strada Academiei 5,
  Chişinău, MD-2028], May 1979.
\newblock In Russian.

\bibitem{Danilcenko1979ParametricalExpressibilitykValuedLogic}
Anna~Fedorovna Daniľčenko.
\newblock On parametrical expressibility of the functions of {$k$}-valued
  logic.
\newblock In \emph{Finite algebra and multiple-valued logic ({S}zeged, 1979)},
  volume~28 of \emph{Colloq. Math. Soc. János Bolyai}, pages 147--159.
  North-Holland, Amsterdam-New York, 1981.

\bibitem{deMouraBjoernerZ3EfficientSMTsolver}
Leonardo de~Moura and Nikolaj Bjørner.
\newblock {Z3:} an efficient {SMT} solver.
\newblock In Cartic~R. Ramakrishnan and Jakob Rehof, editors,
  \emph{{I}nternational {C}onference on {T}ools and {A}lgorithms for the
  {C}onstruction and {A}nalysis of {S}ystems ({TACAS} 2008)}, volume 4963 of
  \emph{Lecture Notes in Comput. Sci.}, pages 337--340, Berlin, Heidelberg,
  March 2008. Springer.
\newblock \doi{10.1007/978-3-540-78800-3_24}.

\bibitem{GeigerClosedSystemsOfFunctionsAndPredicates}
David Geiger.
\newblock Closed systems of functions and predicates.
\newblock \emph{Pacific J. Math.}, 27(1):95--100, 1968.
\newblock Online available from
  \url{http://projecteuclid.org/euclid.pjm/1102985564}.

\bibitem{JanovMucnik1959}
Jurij~Ivanovič Janov and Aľbert~Abramovič
  Mučnik.
\newblock \rus{О существовании $k$-значных замкнутых классов, не имеющих конечного базиса}
 [{O}n the existence of {$k$}-valued closed classes having no finite basis].
\newblock \emph{Dokl. Akad. Nauk SSSR}, 127(1):44--46, 1959.

\bibitem{KuznecovCentralisers1979}
Aleksandr~Vladimirovič Kuznecov.
\newblock \rus{О средствах для обнаружения невыводимости или невыразимости}
  [{M}eans for detection of nondeducibility and inexpressibility].
\newblock In \emph{\rus{Логический вывод} [Logical Inference] ({M}os\-cow, 1974)},
  pages 5--33. ``{N}auka'', {M}os\-cow, 1979.
\newblock In Russian.

\bibitem{PoeKal}
Reinhard Pöschel and Lev~Arkaďevič Kalužnin.
\newblock \emph{Funk\-tio\-nen- und {R}e\-la\-tio\-nen\-al\-ge\-bren.
  {Ein {K}a\-pi\-tel der dis\-kre\-ten {M}a\-the\-ma\-tik} [Function
  and relation algebras. A chapter in discrete mathematics]}, volume~15
  of \emph{Ma\-the\-ma\-ti\-sche {M}o\-no\-gra\-phi\-en [Mathematical
  Monographs]}.
\newblock {VEB} {D}eut\-scher {V}er\-lag der {W}is\-sen\-schaf\-ten, Berlin,
  1979.

\bibitem{Z3}
Microsoft Research.
\newblock {Z3} {T}heorem {P}rover, 2020.
\newblock Available from {\url{https://github.com/z3prover/z3}} or from
  {\url{https://rise4fun.com/Z3/}}.

\bibitem{SnowGeneratingPrimitivePositiveClones}
John~W. Snow.
\newblock Generating primitive positive clones.
\newblock \emph{Algebra Universalis}, 44(1-2):169--185, 2000.
\newblock \doi{10.1007/s000120050179}.

\end{thebibliography}
